\documentclass[12pt,a4paper]{amsart}
\usepackage{amssymb,latexsym,graphicx}
\usepackage{color}
\usepackage[cp1252]{inputenc} %% For Windows
\usepackage[T1]{fontenc}
\usepackage{lmodern}
\usepackage[english]{babel}
\usepackage[section]{placeins}
\usepackage[pagebackref=true]{hyperref}
\hypersetup{pdftex,colorlinks=true,linkcolor=blue,citecolor=blue,urlcolor=blue}
\usepackage{currfile}

\newtheorem{theorem}{Theorem}[section]
\newtheorem{lemma}[theorem]{Lemma}
\newtheorem{proposition}[theorem]{Proposition}
\newtheorem{definition}[theorem]{Definition}

\newtheorem{remark}[theorem]{Remark}
\newtheorem{remarks}[theorem]{Remarks}

\newtheorem{assump}[theorem]{Assumption}
\newtheorem{notation}[theorem]{Notation}

\numberwithin{equation}{section}
\numberwithin{figure}{section}
\numberwithin{table}{section}

\definecolor{purple}{RGB}{127,0,255}

\newcommand{\E}{\mathbb{E}}
\newcommand{\N}{\mathbb{N}}

\newcommand{\R}{\mathbb{R}}

\newcommand{\bS}{\mathbb{S}}
\newcommand{\T}{\mathbb{T}}

\newcommand{\cE}{\mathcal{E}}
\newcommand{\cG}{\mathcal{G}}

\newcommand{\cL}{\mathcal{L}}

\newcommand{\cS}{\mathcal{S}}
\newcommand{\cT}{\mathcal{T}}
\newcommand{\cZ}{\mathcal{Z}}

\newcommand{\mf}{\mathfrak}

\newcommand{\noi}{\noindent}

\newcommand{\ecp}{\mathrm{\textsc{ECP}}}

\newcommand{\wb}[1]{\overline{#1}}
\newcommand{\sm}{\!\setminus\!}

\newcommand{\rpn}[1]{\mathbb{R}\mathrm{P}^{#1}}

\begin{document}

\title[Deformation of level sets. Application.]{ Level sets of certain Neumann eigenfunctions under deformation of Lipschitz domains \\ Application to the Extended Courant Property}

\author[P. B\'{e}rard]{Pierre B\'erard}
\author[B. Helffer]{Bernard Helffer}

\address{PB: Universit\'{e} Grenoble Alpes and CNRS\\
Institut Fourier, CS 40700\\ 38058 Grenoble cedex 9, France.}
\email{pierrehberard@gmail.com}

\address{BH: Laboratoire Jean Leray, Universit\'{e} de Nantes and CNRS\\
F44322 Nantes Cedex, France and LMO (Universit\'e Paris-Sud).}
\email{Bernard.Helffer@univ-nantes.fr}



\begin{abstract}
In this paper, we prove that the Extended Courant Property fails to be true for certain smooth,  strictly convex domains with Neumann boundary condition: there exists a linear combination of a second and a first Neumann eigenfunctions, with   three nodal domains. For the proof, we revisit a deformation argument of Jerison and Nadirashvili  (J. Amer. Math. Soc. 2000, vol. 13). This argument being interesting in itself, we give full details. In particular, we carefully control the dependence of the constants on the geometry of our Lipschitz domains along the deformations.\\[5pt]

\noi \textsc{R\'{e}sum\'{e}.} Dans cet article, nous montrons que la ``propri\'{e}t\'{e} \'{e}tendue de Courant'' est fausse pour certains domaines convexes lisses avec condition au bord de Neumann~: il existe une combinaison lin\'{e}aire d'une premi\`{e}re et d'une seconde fonctions propres de Neumann ayant trois domaines nodaux. Pour la d\'{e}monstration, nous reformulons un argument de Jerison et Nadirashvili  (J. Amer. Math. Soc. 2000, vol. 13). Cet argument \'{e}tant int\'{e}ressant en lui-m\^{e}me, nous d\'{e}taillons la preuve. En particulier, nous explicitons la d\'{e}pendance des constantes par rapport \`{a} la g\'{e}om\'{e}trie des domaines lipschitziens le long des d\'{e}formations.
\end{abstract}%

\date{\today ~(\currfilename)}
%\date{\today}
\maketitle

%\newpage
%\tableofcontents
%\newpage

\section{Introduction}\label{S-intro}

Let $\Omega \subset \R^d$ be a bounded domain (open and connected), with $d\ge 2$. We assume that $\Omega$ is smooth enough, and we consider the eigenvalue problem
\begin{equation}\label{E-intro-2}
\left\{
\begin{array}{l}
- \Delta \varphi = \mu \, \varphi \text{~in~} \Omega\,,\\[5pt]
B(\varphi) = 0 \text{~on~} \partial \Omega\,,
\end{array}
\right.
\end{equation}
where the boundary condition $B(\varphi)$ is either the Dirichlet boundary condition $\varphi|_{\partial \Omega}=0$, or the Neumann boundary condition  $\frac{\partial \varphi}{\partial n_e}|_{\partial \Omega}=0$ (here $n_e$ denotes the exterior unit normal).\medskip

We write the eigenvalues of \eqref{E-intro-2} in nondecreasing order, with multiplicities, starting with the index $1$,
%i
\begin{equation}\label{E-intro-4}
\mu_1(\Omega,\mf{a}) < \mu_2(\Omega,\mf{a}) \le \mu_3(\Omega,\mf{a}) \le \cdots \,,
\end{equation}
where $\mf{a} \in \{\mf{d,n}\}$ denotes the boundary condition. \medskip

Given an eigenvalue $ \mu(\Omega,\mf{a})$ of \eqref{E-intro-2}, we denote by $\cE\left(\mu(\Omega,\mf{a})\right)$ the corresponding eigenspace. Given an eigenfunction $\varphi \in \cE\left(\mu(\Omega,\mf{a})\right)$, we denote by
\begin{equation}\label{E-intro-6}
\cZ(\varphi) = \overline{\{x\in \Omega ~|~ \varphi(x)=0\}}
\end{equation}
the \emph{nodal set} of $\varphi$, and by $\beta_0(\varphi)$ the number of \emph{nodal domains} (the connected components of $\Omega\! \setminus \! \cZ(\varphi)$) of the function $\varphi$.\medskip

Given an eigenvalue $ \mu = \mu(\Omega,\mf{a})$ of \eqref{E-intro-2}, we denote by $ \kappa(\mu)$ the \emph{least index} of $\mu$,
\begin{equation}\label{E-intro-8}
\kappa(\mu) = \min \{k ~|~ \mu_k(\Omega,\mf{a}) = \mu\}\,.
\end{equation}

The following classical theorem was proved by R.~Courant in 1923, see for example \cite[$\S$~VI.6]{CH1953}.

\begin{theorem}[Courant's nodal domain theorem]\label{T-intro-2}
Let $\mu$ be an eigenvalue of \eqref{E-intro-2}, and $\varphi \in \cE(\mu)$ a cor\-res\-ponding eigenfunction. Then,
\begin{equation}\label{E-intro-10}
\beta_0(\varphi) \le \kappa(\mu)\,.
\end{equation}
\end{theorem}%

When $d=1$, given a finite interval $]\alpha,\beta[$,  instead of the eigenvalue problem for the Laplacian, we consider the Sturm-Liouville eigenvalue problem,
\begin{equation}\label{E-intro-20}
\left\{
\begin{array}{l}
- y'' + q\, y = \mu \, y \text{~in~} ]\alpha,\beta[\,,\\[5pt]
B(y) = 0 \text{~at~} \{\alpha,\beta\}\,,
\end{array}%
\right.
\end{equation}
where $q$ is a smooth real function on $[\alpha,\beta]$. There are striking differences between the eigenvalue problems~\eqref{E-intro-20} ($d=1$) and  \eqref{E-intro-2} ($d \ge 2$).\medskip

\noi \emph{First difference.} \\
When $d=1$, a classical theorem of C.~Sturm \cite{Sturm1836a} states that the eigenvalues of \eqref{E-intro-20} are all simple, and that an eigenfunction of \eqref{E-intro-20}, associated with the $n$th eigenvalue, has exactly $n$ nodal domains. \medskip

When $d\ge 2$, the eigenvalues of \eqref{E-intro-2} may have multiplicities (this is for example the case for a square with either Dirichlet or Neumann condition on the boundary). By Courant's nodal domain theorem, an eigenfunction of \eqref{E-intro-2}, associated with the $n$th-eigenvalue has at most $n$ nodal domains. However,
\begin{enumerate}
  \item  For the round sphere $\bS^2$, and for the square with Dirichlet boundary condition, examples of A.~Stern \cite{BH-tsg1,BH-tsg2} show that  there is no general lower bound on $\beta_0(\varphi)$ for higher energy eigenfunctions, except the trivial bound $\beta_0(\varphi) \ge 2$\,.  Note that  the example of the square suggests that such a statement might not be true for the Neumann boundary condition,  see the paragraph before Proposition~10.2 in \cite{HePe2014}.
  \item A theorem of {\AA}.~Pleijel \cite{Ple1956} shows that the upper bound $\beta_0(\varphi) \le \kappa(\mu)$ is sharp for finitely many eigenvalues $\mu$ only.
\end{enumerate}

\noi \emph{Second difference.}\\
Another, not so well-known,  theorem of C.~Sturm \cite{Sturm1836b}  states that, for $n \ge m \ge 1$, a linear combination $\sum_{k=m}^n a_k V_k$ of eigenfunctions of \eqref{E-intro-20}, in the range $k \in \{m,\ldots,n\}$, has at least $(m-1)$, and at most $(n-1)$ zeros in the interval $]\alpha,\beta[$.  We refer to \cite{BH-sturm} for a more precise statement of Sturm's theorem,  and to \cite{GaKr2002}, in particular Theorem~1 in Section~IV.3, for a different point of view.\smallskip

In dimension $d\ge 2$, a similar statement (for the upper bound) appears in Footnote~1, page 454 of \cite[Chap.~VI.6]{CH1953}, namely:\\ \emph{Any linear combination of the first $n$ eigenfunctions divides the domain, by means of its nodes, into no more than $n$ subdomains. See the G\"{o}ttingen dissertation of H.~Herrmann, Beitr\"{a}ge zur Theorie der Eigenwerte und Eigenfunktionen, 1932.} \\
This statement is sometimes referred to as the ``Courant-Herrmann theorem'' \cite[$\S$~9.2]{Gud1974}, or the ``Courant-Herrmann conjecture'' \cite{GlZh2003}. We shall call this statement the \emph{Extended Courant Property}, and refer to it as the $\ecp(\Omega,\mf{a})$, when applied to the boundary value problem \eqref{E-intro-2}, with the boundary condition $\mf{a}$. \medskip

 In \cite{Arn2014},  see also \cite{Arn2011,Kuz2015},  V.~Arnold points out that the $\ecp(\rpn{2},g_0)$ is true for the round metric $g_0$, and that the $\ecp(\rpn{3},g_0)$ is false, with counterexamples constructed by  O.~Viro \cite{Vir1979}. Arnold also claims that $\ecp(\bS^2,g)$ is false for a generic metric $g$.  As far as we understand, the only  published proof that  the assertion ``the $\ecp(\rpn{2},g_0)$ is true'', is  the real algebraic geometry proof given in \cite{Ley1996} (Theorem~1, and second remark on page 305).  To our knowledge, no proof  of the second claim has been published, see \cite[Section~5]{BCH2019} for a related result. \medskip

Little seems to be known on the $\ecp$. In \cite{BH-ecp1,BH-ecp2}, we  give some examples of domains such that $\ecp(\Omega,\mf{a})$ is false, with either the Dirichlet or the Neumann boundary condition. However, all these examples are singular (domains or surfaces with  cracks), or have a nonsmooth boundary (polygonal domains). A natural question is whether one can construct counterexamples to the $\ecp$ with  a $C^{\infty}$ boundary. Numerical simulations for the equilateral triangle with rounded corners (the corners of the triangle are replaced by circular caps tangent to the sides) suggest that this should be true. Note however that a triangle with rounded corners is $C^1$, not $C^2$.\medskip

\begin{figure}[!ht]
  \centering
  \includegraphics[scale=0.45]{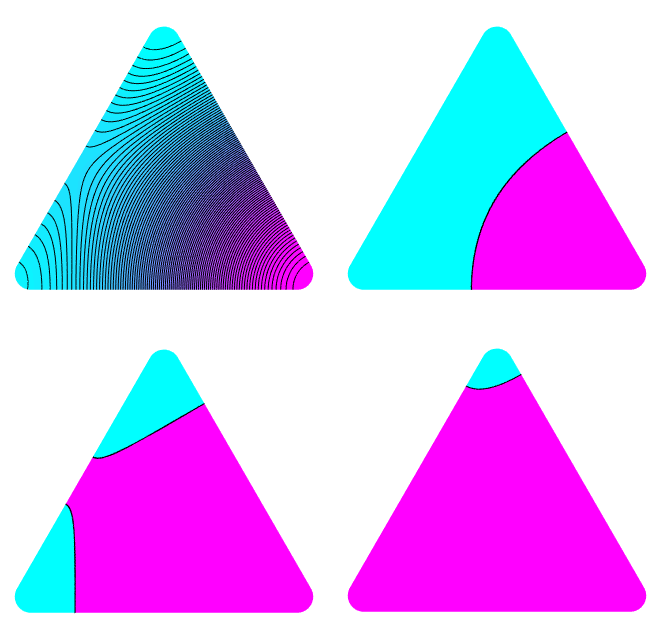}
\caption{Level sets of one of the second  Neumann eigenfunctions of the equilateral triangle with rounded corners}\label{RET}
\end{figure}

 The pictures in the first row of Figure~\ref{RET} display the level sets and nodal domains of a second Neumann eigenfunction $\phi$ of the equilateral triangle with rounded corners, as calculated by \textsc{matlab}. The function is almost symmetric\footnote{Generally speaking, numerical softwares do not necessarily produce the symmetric eigenfunctions when an eigenvalue is not simple.} with respect to one of the axes of symmetry of the triangle. The pictures in the second row display the nodal sets of the function $a + \phi$ for two values of $a$. They provide a numerical evidence that $\ecp$ is not true for the equilateral triangle with rounded corners, and Neumann boundary condition. \medskip

\FloatBarrier

In this paper, we prove,

\begin{theorem}\label{T-intro-4}
There exists a one-parameter family of $C^{\infty}$,  strictly convex domains $\{\Omega_t, 0 < t < t_0\}$ in $\R^2$, with the symmetry of the equilateral triangle $\cT_e$, such that:
\begin{enumerate}
  \item The family is strictly  increasing, and $\Omega_t$ tends to $ \cT_e$,  in the sense of the Hausdorff distance, as $t$ tends to $0$.
  \item For any $t \in ]0,t_0[$,  the  $\ecp(\Omega_{t},\mf{n})$ is false. More precisely, for each $t$, there exists a linear combination of a symmetric 2nd Neumann eigenfunction and a 1st Neumann eigenfunction of $\Omega_{t}$, with  precisely three nodal domains.
\end{enumerate}
\end{theorem}%

The starting point of the proof of Theorem~\ref{T-intro-4} is the fact, established in \cite{BH-ecp1}, that the $\ecp(\cT_e,\mf{a})$ is false for both the Dirichlet, and the Neumann boundary conditions on the equilateral triangle. The idea  is then to  show that one can find a deformation of $\cT_e$  by smooth strictly convex domains, in such a way that the symmetric second Neumann eigenfunction deforms nicely.

\subsection*{Organization of the paper}
 In Section~\ref{S-jnn},  we revisit a deformation argument given by Jerison and Nadirashsvili \cite{JeNa2000} in the framework of the ``hot spots'' conjecture. The main result is Lemma~\ref{L-jnn-32s}. This argument being interesting in itself, we give full details. In Section~\ref{S-teqal}, we construct smooth strictly convex approximations of the equilateral triangle by using the convexity properties of its first Dirichlet eigenfunction (Proposition~\ref{P-teqa-0d}) or its torsion function (Proposition~\ref{P-teqa-0f}). These approximating domains have the symmetries of the equilateral triangle. A key point is that their second Neumann eigenspace has dimension $2$, with a nice symmetry property (Proposition~\ref{P-teqa-n2}).
In Section~\ref{S-appl}, we first construct yet another deformation $\{\Omega_t\}$ of the equilateral triangle (Proposition~\ref{P-appln-2}), and then complete the proof of Theorem~\ref{T-intro-4} using this deformation. %

\begin{remark}\label{R-add}
As pointed out by the anonymous referee, a natural question arises from the counterexamples to the Extended Courant Property. \emph{Does there exist a constant $C$ such that every linear combination of the first $n$ eigenfunctions has at most $C\, n$ nodal domains, for some constant $C$.} The answer is \textsc{no} in general. The first examples are constructed in \cite{BH-ecp1}, Remark~4.3 and 6.2, by introducing cracks. Further examples are constructed in \cite{BLS2018} for the $2$-torus, and in \cite{BCH2019} for regular polygons with Neumann boundary condition, for the $2$-torus and for the $2$-sphere. More precisely, in the latter case, there exists a metric $g$ on $\T^2$ (resp. $\bS^2$), and an associated eigenfunction $\Phi$ of the Laplace-Beltrami operator $\Delta_g$, such that the set $\{\Phi > 1\}$ has infinitely many connected components. Furthermore, the metric $g$ can be chosen as close as desired from the flat (resp. the round) metric.
\end{remark}%

\subsection*{Acknowledgements} The authors would like to thank P.~Bousquet,  T.~Hoffmann-Ostenhof, and H.~Tamura,  for providing useful references.

\section{ A deformation argument}\label{S-jnn}

In this section, we revisit a deformation argument of Jerison and Nadirashvili \cite[Section~2]{JeNa2000}.  Note that our framework is different: they are interested in antisymmetric eigenfunctions in domains with two orthogonal lines of mirror symmetry; we are interested in symmetric eigenfunctions in domains with the symmetries of an equilateral triangle.  Because we work with symmetric eigenfunctions, we need an extra assumption (Assumption~\ref{A-jnn-9s}) which is satisfied by the domains used in the proof of Theorem~\ref{T-intro-4}, Section~\ref{S-appl}. We also aim at controlling the constants which appear in the analytic inequalities, and at making sure that they are uniform in a large class of domains. This aspect is not always taken care of clearly in the literature.\medskip

\subsection{Geometric framework: the class $\cL_M$}\label{SS-jnn-1}

Let $M$ be a positive constant.

\begin{definition}\label{A-jnn-LM} The class $\cL_M$ comprises the sets $\Omega \subset \R^2$ which satisfy the following conditions.
\begin{equation}\label{E-jnn-2}
\Omega \text{~is convex and open, with~} 0 \in \Omega\,.
\end{equation}
\begin{equation}\label{E-jnn-3}
\wb{B}(M^{-1}) \subset \Omega \subset \wb{\Omega} \subset B(M)\,,
\end{equation}
where $B(R)$ denotes the open ball centered at $0$, with radius $R$, and $\wb{B}(R)$ denotes the corresponding closed ball.
\begin{equation}\label{E-jnn-4}
\Omega \text{~is symmetric with respect to~} D:= \{(u,v) \in \R^2 ~|~ u=0\}\,.
\end{equation}
\begin{equation}\label{E-jnn-5}
\partial \Omega \text{~is regular at~} D \cap \partial \Omega\,,
\end{equation}
i.e. in a neighborhood of $m \in D\cap \partial \Omega$, the boundary $\partial \Omega$ is piecewise $C^1$, and $\partial \Omega \sm\{m\}$ is $C^1$.\smallskip

The domain $\Omega$ can be described by a polar equation,
\begin{equation}\label{E-jnn-6}
\Omega = \left\{  (r,\theta) ~|~ 0 \le r < \rho(\theta)\right\}\,,
\end{equation}
where the function $\rho$ is a $2\pi$-periodic, Lipschitz function, with Lipschitz constant bounded from above by $M$.
\end{definition}%

We define the domain,
\begin{equation}\label{E-jnn-7b}
\Omega_{+} := \Omega \cap \{(u,v) \in \R^2 ~|~ u > 0\}\,.
\end{equation}
We decompose its boundary $\partial \Omega_{+}$ as
\begin{equation}\label{E-jnn-7c}
\partial \Omega_{+} = \overline{\Gamma \sqcup \Gamma_D}\,,
\end{equation}
with $\Gamma = \partial \Omega_{+}\cap \{u > 0\}$, and $\Gamma_D = D \cap \Omega_{+}$.\medskip

\textbf{Notation.} In the sequel, we denote by $D$ both the line, and the mirror symmetry with respect to the line $D$. We denote by $D^*$ the action of the symmetry $D$ on functions, $D^*\phi = \phi\circ D$.

\begin{remarks}\label{R-jnn1-2} We note the following properties for later reference.
\begin{enumerate}
  \item According to Proposition~2.4.4 in \cite{HePi2005}, domains satisfying conditions~\eqref{E-jnn-2} and \eqref{E-jnn-3} satisfy a uniform (i.e. depending only on $M$) cone property. It follows from Theorem~2.4.7, and Remark~2.4.8 in \cite{HePi2005} that such domains are uniformly Lipschitz domains (i.e., the boundary is locally the graph of a Lipschitz function, ibidem Definition~2.4.5).
  \item With the definitions of \cite{Gri1985}, for such domains, the inclusion $H^1(\Omega) \hookrightarrow L_2(\Omega)$ is compact, and we can define eigenvalues using the variational approach.
  %% See also Appendix~\ref{A-be}.
  \item The fact that a domain $\Omega$, defined in polar coordinates as in \eqref{E-jnn-6}, is a Lipschitz domain also follows from \cite[Theorem~7.1]{Wue2008}.
  \item Let $\Omega$ be a domain defined by a polar equation, as in \eqref{E-jnn-6}. Define the function $r(\theta)$ by $r(\theta) = 1/\rho(\theta)$. If $\Omega$ is convex, then the second derivative of $r$, in the sense of distributions, is a measure such that $r''(\theta) + r(\theta) \ge 0\,$, see \cite[Chap.~3.4]{Egg1958}. %%See Appendix~\ref{A-conv}.
\end{enumerate}
\end{remarks}%

We consider the Neumann eigenvalue problem for $-\Delta$ in $\Omega$. We denote the Neumann eigenvalues by $\nu_i(\Omega)$, and arrange them in nondecreasing order, starting with the index $1$. We also consider the eigenvalue problems for $-\Delta$ in $\Omega_{+}$, with either the Neumann boundary condition on $\partial \Omega_{+}$, or the mixed boundary conditions, Neumann on $\Gamma$ and Dirichlet on $\Gamma_D$. We denote these eigenvalues respectively by $\mu_i(\Omega_{+},\mf{nn})$, and $\mu_i(\Omega_{+},\mf{nd})$, and arrange them in nondecreasing order, starting with the index $1$.
\medskip

We are interested in the \emph{least positive} eigenvalues of $\Omega$ associated with the symmetry $D$. More precisely, we introduce
\begin{equation}\label{E-jnn-9a}
\nu^{-}(\Omega) := \inf \{\nu_i(\Omega) ~|~ i \ge 2\,, ~\exists \varphi\,, -\Delta \varphi = \nu_i(\Omega) \varphi\,, ~D^*\varphi = - \varphi\}\,,
\end{equation}
and
\begin{equation}\label{E-jnn-9s}
\nu^{+}(\Omega) := \inf \{\nu_i(\Omega) ~|~ i \ge 2\,, ~\exists \varphi\,, -\Delta \varphi = \nu_i(\Omega) \varphi\,, ~D^*\varphi = \varphi\}\,,
\end{equation}
where the equations $-\Delta \varphi = \nu_i(\Omega) \varphi$ are to be understood in $\Omega$.\medskip

It is easy to see that
\begin{equation}\label{E-jnn-9c}
\left\{
\begin{array}{lll}
\nu^{-}(\Omega) &=& \mu_1(\Omega_{+},\mf{nd})\,,\\[5pt]
\nu^{+}(\Omega) &=& \mu_2(\Omega_{+},\mf{nn})\,,\\[5pt]
\nu_2(\Omega) &=& \min \{\nu^{-}(\Omega) \,, \nu^{+}(\Omega) \}\,.
\end{array}
\right.
\end{equation}

\begin{remarks}\label{R-jnn-9}About the eigenvalues $\nu^{-}(\Omega) $ and $\nu^{+}(\Omega)$.
\begin{enumerate}
   \item  Because $\mu_1(\Omega_{+},\mf{nd})$ is simple,  there is, up to scaling\footnote{By this, we mean ``up to multiplication by a nonzero scalar''.}, a unique anti-symmetric eigenfunction $\psi_{\Omega}$ of $-\Delta$ in $\Omega$, associated with the eigenvalue $\nu^{-}(\Omega)$,
\begin{equation}\label{E-jnn-9d}
\int_{\Omega} \psi^2_{\Omega} = 1 \text{~and~} \psi_{\Omega}|_{\Omega_{+}} > 0\,.
\end{equation}
  \item If $\nu_2(\Omega)$ is a simple eigenvalue, then either $\nu_2(\Omega) = \nu^{+}(\Omega) < \nu^{-}(\Omega)$ or ~$\nu_2(\Omega) = \nu^{-}(\Omega) < \nu^{+}(\Omega)$, and the corresponding eigenfunction is either invariant, or anti-invariant under $D$.
   \item If $\dim \cE\big( \nu_2(\Omega) \big) \ge 2$, then
    $$
    \cE(\nu_2) = \left(\cE(\nu_2) \cap \cS_{+}\right)  \bigoplus  \left(\cE(\nu_2) \cap \cS_{-}\right),
$$
with $ \dim \cE(\nu_2) \cap \cS_{-} \le 1$.  Here, we have used the notation
\begin{equation}\label{E-jnn-9not}
\cS_{\sigma} := \{\phi ~|~ D^{*}\phi = \sigma\, \phi\}\,, ~\sigma \in \{+,-\}\,.
\end{equation}
\item  If $\Omega$ is sufficiently regular, then $\dim \cE(\nu_2) \le 3$, see \cite{Che1976,HoMiNa1999}.
  \item Let $\cT_i(\alpha)$ be an isosceles triangle with aperture $\alpha \in ]0,\pi[$. According to \cite[$\S$~10]{LaSi2017},
$$
\nu_2(\cT_i(\alpha)) = \nu^{+}(\cT_i(\alpha)) < \nu^{-}(\cT_i(\alpha)) \text{~when~}  0 < \alpha < \frac{\pi}{3}\,,
$$
$$
\nu_2(\cT_i(\alpha)) = \nu^{-}(\cT_i(\alpha)) < \nu^{+}(\cT_i(\alpha)) \text{~when~}  \frac{\pi}{3} < \alpha < \pi \,.
$$
There is a bifurcation at $\frac{\pi}{3}$, in which case
$$
\nu_2(\cT_i(\frac{\pi}{3}) )= \nu^{-}(\cT_i(\frac{\pi}{3})) = \nu^{+}(\cT_i(\frac{\pi}{3}))\,.
$$
%\item Numerical computations for the equilateral rhombus suggest that
%$\nu^{+}$ may have multiplicity $2$, and that $\nu^{+}$ and $\nu^{-}$
%are not necessarily consecutive eigenvalues.
 \item  In Section~\ref{S-teqal}, we consider domains $\Omega$ which admit the symmetry group $\cG_0$ of the equilateral triangle, see \eqref{E-teqa-2}. For such domains, Proposition~\ref{P-teqa-n2} tells us that
$$
 \nu^{-}(\Omega) = \nu^{+}(\Omega) = \nu_2(\Omega) = \nu_3(\Omega) < \nu_4(\Omega).
$$
\end{enumerate}\vspace{-3mm}
\end{remarks}%

\textbf{Notation.}~In \eqref{E-jnn-9d}, and henceforth, we skip the (Lebesgue) measure $dx$ in the integrals.\medskip

 We now introduce a technical assumption.

\begin{assump}\label{A-jnn-9s}
The eigenvalue $\mu_2(\Omega_{+},\mf{nn})$ is simple.
\end{assump}%

 Remark~\ref{R-jnn-9}-(6) tells us that Assumption~\ref{A-jnn-9s}  is satisfied by convex domains with the $\cG_0$ symmetry, see Proposition~\ref{P-teqa-n2}, and in particular by the domains $\Omega_t$ constructed for the proof of Theorem~\ref{T-intro-4}.

\begin{remark}\label{R-jnn-9s}
Provided that Assumption~\ref{A-jnn-9s} is satisfied, there is a $D$-symmetric  eigenfunction $\phi_{\Omega}$ of $-\Delta$ in $\Omega$, associated with $\nu^{+}(\Omega)$. This eigenfunction is uniquely determined, up-to-sign, by the normalization $\int_{\Omega} \phi^2_{\Omega} = 1$.  In Lemma~\ref{L-jnn-32s}, we will prove that one can actually make a unique choice of $\phi_{\Omega_t}$ along a path of domains.
\end{remark}%

%%\^{A}£

\subsection{Preliminary estimates}\label{SS-jnn-2}

We shall now examine how the eigenvalues $\nu^{\pm}(\Omega)$, and the corresponding eigenfunctions, vary with the domain $\Omega \in \cL_M$. For this purpose, and following \cite{JeNa2000}, we introduce the following distance in the class $\cL_M$,
\begin{equation}\label{E-jnn-10}
d_r(\Omega_1,\Omega_2) = \|\rho_1-\rho_2\|_{\infty}\,,
\end{equation}
if the domains are defined by the functions $\rho_1$ and $\rho_2$ respectively, as in \eqref{E-jnn-6}.\medskip

Note that this distance is bigger than the Hausdorff distance between open sets contained in a given compact ball $D$,
\begin{equation}\label{E-jnn-10a}
d_H(\Omega_1,\Omega_2) := d^H(D\sm\Omega_1,D\sm\Omega_2) \,.
\end{equation}
Here,
\begin{equation}\label{E-jnn-10aa}
d^H(K_1,K_2) := \max \left\lbrace \sup_{x\in K_1} \inf_{y\in K_2} d(x,y) \,, \sup_{x\in K_2} \inf_{y\in K_1} d(x,y)\right\rbrace \,,
\end{equation}
is the Hausdorff distance between the compact sets $K_1$ and $K_2$, and
$d(x,y)$ is the Euclidean distance between the points $x,y \in \R^2$.\medskip

Note that the distance defined in \eqref{E-jnn-10a} does not depend on the choice of the compact $D$, once it contains both $\Omega_1$ and $\Omega_2$.\medskip

 \noi \textbf{Notation}. In the sequel, $|\Omega|$ denotes the area of a domain $\Omega$. We will also use the following convention. We use constants $C_i, i \in \N$ in the statements, and local constants $C_{i,j}, i,j \in \N$   inside  the proofs. Note that the constants are not numbered linearly. When a constant appears, we mention which parameters it depends upon.

\begin{lemma}\label{L-jnn-2}
There exists a constant $C_1(M)$ such that, for any domains $\Omega_1, \Omega_2 \in \cL_M$,
\begin{equation}\label{E-jnn-L2}
|\Omega_1 \setminus \Omega_2| \le C_1(M)\, d_r(\Omega_1,\Omega_2)\,.
\end{equation}
\end{lemma}%

\proof It suffices to notice that
\begin{equation*}
\Omega_1 \setminus \Omega_2 = \{(r,\theta) ~|~ \rho_2(\theta) \le r < \rho_1(\theta)\}\,,
\end{equation*}
and to compute the area in polar coordinates.\hfill \qed

\begin{lemma}\label{L-jnn-4}
There exists a constant $C_2(M)$ such that, for any $\Omega \in \cL_M$,
\begin{equation}\label{E-jnn-L4}
\max \{ \nu_2(\Omega)\,, \nu^{+}(\Omega)\,, \nu^{-}(\Omega) \} \le C_2(M)\,.
\end{equation}
\end{lemma}%

\proof  Since $\Omega \in \cL_M$, condition \eqref{E-jnn-3} is satisfied. We then have,
\begin{equation*}
\left\{
\begin{array}{l}
\nu_2(\Omega) \le \delta_2(\Omega) < \delta_2(B(M^{-1}))\,,\\[5pt]
\nu^{+}(\Omega) = \mu_2(\Omega_{+},\mf{nn}) \le \delta_2(\Omega_{+})
\le \delta_2\big( B(M^{-1})\cap \{u > 0\}\big) \,,\\[5pt]
\nu^{-}(\Omega) = \mu_1(\Omega_{+},\mf{nd}) \le \delta_1(\Omega_{+})
\le \delta_1\big( B(M^{-1})\cap \{u > 0\}\big) \,,
\end{array}
\right.
\end{equation*}
where we have used $\delta$'s to denote Dirichlet eigenvalues. \hfill \qed

\begin{proposition}\label{P-jnn-6}
 Under the Assumption~\ref{A-jnn-9s}, there exists a constant $C_3(M)$ such that, for any $\Omega \in \cL_M$, the normalized eigenfunction $\psi_{\Omega}$ (defined in Remark~\ref{R-jnn-9}-(1)), and the normalized eigenfunction $\phi_{\Omega}$ (defined in Remark~\ref{R-jnn-9s}), belong to the Sobolev space $H^{2}(\Omega)$, with corresponding Sobolev norm less than or equal to $C_3(M)$,
\begin{equation}\label{E-jnn-P6}
\|\psi_{\Omega} \|_{H^{2}(\Omega)}  + \|\phi_{\Omega} \|_{H^{2}(\Omega)} \le C_3(M)\,.
\end{equation}
\end{proposition}%

\proof  We refer to \cite{Gri1985}, proofs of Theorem 3.2.1.2 and 3.2.1.3.  The point we want to stress here, is that the bound is uniform with respect to the domains in $\cL_M$.
%% We sketch the proof in Appendix~\ref{A-be}.
\hfill \qed

\begin{remark}\label{R-be-R-jena}
The $H^2$ estimates in the proposition hold for \emph{convex} domains. For more general Lipschitz domains, there are only $H^s$ estimates, with $s=\frac{3}{2}$ in \cite{JeNa2000}, or $s < \frac{3}{2}$ in \cite{Sav1998}. A counterexample is given in \cite{Gri1985}.
\end{remark}%

\begin{proposition}[Extension theorem]\label{P-jnn-8}
For any domain $\Omega \in \cL_M$, there exists a linear extension operator $\E_{\Omega}$, such that for any $s > 0$,
\begin{equation*}
\E_{\Omega} : H^{s}(\Omega) \to H^{s}(\R^n)\,,
\end{equation*}
and there exists a posi\-tive constant $C_4(M,s)$, such that, for all $\varphi \in H^{s}(\Omega)$,
\begin{equation}\label{E-jnn-L8}
\left\{
\begin{array}{l}
\|\E_{\Omega}(\varphi)\|_{H^{s}(\R^n)} \le C_4(M,s) \|\varphi\|_{H^{s}(\Omega)}\,,\\[5pt]
\E_{\Omega}(\varphi)|_{\Omega} = \varphi \text{~~almost everywhere}\,,\\[5pt]
\E_{\Omega}(\varphi) \text{~is~} D\text{-(anti)symmetric, if~} \varphi \text{~is.}
\end{array}
\right.
\end{equation}
Furthermore, one can choose $\E_{\Omega}(\varphi)$ with compact support in $B(2M)$.
\end{proposition}%

\proof This proposition follows from Theorem~5 in \cite[Chap.~VI.3]{Ste1970} and interpolation.  We again point out that the constant $C_4(M,s)$ is uniform in $\cL_M$. \hfill \qed \medskip

Finally, we mention the classical Sobolev embedding theorem, in the form  we will use later on. Recall that $B(R)$ is the open ball with center the origin, and radius $R$ in $\R^2$.

\begin{proposition}\label{P-jnn-10}
For all $\alpha \in [0,1[$, the space $H^2(B(R))$ embeds continuously in $C^{0,\alpha}(\wb{B}(R))$. The space $H^1(B(R))$ embeds continuously in $L_p(B(R))$ for all $p\ge 2$. In particular, for any $s, ~1 \le s < 2$, and for any $\varphi \in H^2(B(R))$, we have $\varphi \in C^{0,s-1}(\wb{B}(R))$, $d\varphi \in L_{\frac{2}{2-s}}(B(R),\R^2)$, and there exists a constant $C_5(R,s)$, such that
\begin{equation}\label{E-jnn-P10}
\|\varphi\|_{L_{\infty}(B(R))} + \|d\varphi\|_{L_{\frac{2}{2-s}}(B(R))} \le C_5(R,s) \|\varphi\|_{H^{2}(B(R))} \,.
\end{equation}
\end{proposition}%

\proof See \cite{Gri1985}, Theorem~1.4.4.1, and equations (1,4,4,3)--(1,4,4,6), for the statements, and Adams \cite{Ada1975}, Chap.~{IV} and {V}, for the proofs. \hfill \qed

\begin{notation}\label{A-jnn-not}
From now on, we choose some $s_0 \in ]1,2[$, and use the notation,
$$
p_0 := p(s_0) = \frac{2}{2-s_0}\,, \text{~and~~} q_0 := q(s_0) = s_0-1 > 0\,.
$$
\end{notation}%

\subsection{Properties of $\nu^{+}(\Omega)$ and $\phi_{\Omega}$}\label{SS-jnn-4}

In this section, we are interested in how the $D$-symmetric eigenfunction $\phi_{\Omega}$ changes along a deformation $\Omega_t$ of the domain.
Note that in \cite{JeNa2000}, Jerison and Nadirashvili consider the $D$-anti-invariant eigenfunctions,  in the context of the ``hot spots'' conjecture.

\begin{lemma}\label{L-jnn-22s}
There exists a constant $C_{20}(M,s_0)$ such that, for any domains $\Omega_1, \Omega_2 \in \cL_M$,
\begin{equation}\label{E-jnn-L22s}
\big| \nu^{+}(\Omega_1) - \nu^+(\Omega_2) \big| \le C_{20} \, d_r(\Omega_1,\Omega_2)^{q_0}\,.
\end{equation}
\end{lemma}%

\proof For the proof, we use the following notation: $\lambda_i = \nu^{+}(\Omega_i)$; $\phi_i = \phi_{\Omega_i}$ is a normalized $D$-invariant eigenfunction of $-\Delta$ in $\Omega_i$, belonging to $\nu^{+}(\Omega_i)$, in particular we have $\int_{\Omega_i} \phi_i = 0$;
$\Phi_i = \E_{\Omega_i}(\phi_{\Omega_i})$ is a $D$-invariant extension of $\phi_{\Omega_i}$, given by Proposition~\ref{P-jnn-8}. We also introduce the function $\Theta_2$ such that
\begin{equation}\label{E-jnn-22sa}
\Theta_2 = \Phi_2 - |\Omega_1|^{-1}\, \int_{\Omega_1} \Phi_2 \,,
\end{equation}
so that $\int_{\Omega_1} \Theta_2 = 0$, and $d\Theta_2 = d\Phi_2$. \medskip

Then,
\begin{equation}\label{E-jnn-22sb}
\int_{\Omega_1} \Theta_2^2 = \int_{\Omega_1} \Phi_2^2 - |\Omega_1|^{-1} \, \left( \int_{\Omega_1} \Phi_2 \right)^2\,.
\end{equation}

Writing
\begin{equation*}
\int_{\Omega_1} \Phi_2 = \int_{\Omega_2}\Phi_2 +  \int_{\Omega_1 \setminus \Omega_2} \Phi_2 -  \int_{\Omega_2 \setminus \Omega_1} \Phi_2 \,,
\end{equation*}
using the fact that $\int_{\Omega_2} \Phi_2 = \int_{\Omega_2} \phi_2 = 0$, Lemma~\ref{L-jnn-2},  Propositions~\ref{P-jnn-6}, \ref{P-jnn-8}, and \ref{P-jnn-10}, we obtain,
\begin{equation*}
\left| \int_{\Omega_1} \Phi_2  \right| \le \|\Phi_2\|_{\infty}\left( |\Omega_1 \setminus \Omega_2| + |\Omega_2 \setminus \Omega_1|  \right)\,,
\end{equation*}
so that there exists a constant $C_{20,1}(M,s_0)$ such that
\begin{equation}\label{E-jnn-22sc}
\left| \int_{\Omega_1} \Phi_2  \right| \le C_{20,1} \, d_r(\Omega_1,\Omega_2) \,.
\end{equation}

We also have
\begin{equation*}
\int_{\Omega_1} \Theta_2^2 = \int_{\Omega_2}\Phi_2^2 +  \int_{\Omega_1 \setminus \Omega_2} \Phi_2^2 -  \int_{\Omega_2 \setminus \Omega_1} \Phi_2^2 - |\Omega_1|^{-1}\, \left( \int_{\Omega_1}\Phi_2  \right)^2 \,.
\end{equation*}
 Using the same arguments as above, as well as \eqref{E-jnn-3}, we obtain that there exists a constant $C_{20,2}(M,s_0)$ such that
\begin{equation}\label{E-jnn-22sd}
1 - C_{20,2} \, d_r(\Omega_1,\Omega_2) \le \int_{\Omega_1} \Theta_2^2 \, \le 1 + C_{20,2} \, d_r(\Omega_1,\Omega_2)\,.
\end{equation}

Similarly, we write
\begin{equation}\label{E-jnn-22sea}
\int_{\Omega_1}|d\Phi_2|^2 = \int_{\Omega_2}|d\Phi_2|^2 + \int_{\Omega_1\setminus \Omega_2}|d\Phi_2|^2 - \int_{\Omega_2\setminus \Omega_1}|d\Phi_2|^2\,.
\end{equation}

Because $(d\Phi_2)|_{\Omega_2} = d\phi_2$, the first integral in the right-hand side is equal to $\lambda_2$. Letting $\Omega$ be either $\Omega_1\setminus \Omega_2$, or $\Omega_2\setminus \Omega_1$, we can write
\begin{equation}\label{E-jnn-22seb}
\int_{\Omega}|d\Phi_2|^2 \le \left( \int_{\Omega}|d\Phi_2|^{2/(2-s_0)} \right)^{2-s_0}\, |\Omega|^{q_0}\,,
\end{equation}
with the Notation~\ref{A-jnn-not}.

As above, recalling that $d\Theta_2 = d\Phi_2$, we conclude that there exists a constant $C_{20,3}(M,s_0)$ such that
\begin{equation}\label{E-jnn-22se}
\int_{\Omega_1} |d\Theta_2|^2 \le \lambda_2 + C_{20,3}\, d_r(\Omega_1,\Omega_2)^{q_0}\,.
\end{equation}

By symmetry between $\lambda_1$ and $\lambda_2$, this completes the proof of the lemma. \hfill \qed\bigskip

 We now consider a family $\{\Omega_t\}_{0 \le t \le a}$ of domains in the class $\cL_M$. We use the notation,
\begin{equation}\label{E-jjn-nota}
\Omega_{t,+} := \Omega_t \cap \{u>0\}\,,
\end{equation}
and we decompose the boundary $\partial \Omega_{t,+}$ into two parts, $\partial \Omega_{t} \cap \{u>0\}$ and $D \cap \Omega_{t,+}$. We assume furthermore that the domains $\Omega_t$ satisfy the Assumption~\ref{A-jnn-9s}, i.e., that the eigenvalues $\nu^{+}(\Omega_t)$, or equivalently  the eigenvalues $\mu_2(\Omega_{t,+},\mf{nn})$, are simple.\medskip

Call $\phi_{t}$ an eigenfunction associated with $\nu^{+}(\Omega_t)$, with $L_2$-norm $1$. It is uniquely defined up to sign. Denote its extension $\E_{\Omega_t}(\phi_t)$ by $\Phi_t$ (see, Proposition~\ref{P-jnn-8}). Recall that $\phi_t$ and $\Phi_t$ are both symmetric with respect to $D$.\medskip

We also use the notation,
\begin{equation}\label{E-jjn-notb}
\left\{
\begin{array}{l}
\lambda_t := \nu^{+}(\Omega_t) = \mu_2(\Omega_{t,+},\mf{nn})\,,\\[5pt]
\mu_0 := \mu_3(\Omega_{0,+},\mf{nn})\,,
\end{array}
\right.
\end{equation}
Observe  that Assumption~\ref{A-jnn-9s} on $\Omega_0$ implies that
\begin{equation}\label{E-jjn-notc}
\lambda_0 < \mu_0 \,.
\end{equation}

 \begin{lemma}\label{L-jnn-32s}
Let $\{\Omega_t\}_{0 \le t \le a}$ be a family of domains in the class $\cL_M$, satisfying  Assumption~\ref{A-jnn-9s}. Assume that $d_r(\Omega_t,\Omega_0)$ tends to zero when $t$ tends to zero.
\begin{enumerate}
  \item For $d_r(\Omega_t,\Omega_0)$ small enough, the function $\phi_t$ can be uniquely defined by the normalization
$$
\int_{\Omega_t} \phi^2_t = 1 \text{~and~} \int_{\Omega_t \cap \Omega_0} \phi_t \phi_0 > 0\,.
$$
\item When $t$ tends to zero, $\Phi_t|_{\Omega_0}$ tends to $\phi_0$ in $L_2(\Omega_0)$. Furthermore, the family $\Phi_t$ is relatively compact in $C^{0,s_0-1}(\R^2)$, and weakly compact in $H^2(\R^2)$.
\item For any $k\in \N$, and for any compact $K \subset \Omega_0$, the functions $\Phi_t$ tend to $\phi_0$ in $C^k(K)$.
\end{enumerate}
\end{lemma}%

\emph{Proof of Assertion (1).}~ We begin as in the proof of Lemma~\ref{L-jnn-22s}. For the time being, $\phi_t$ is well-defined up to sign. Let
\begin{equation}\label{E-jnn-42a}
\Theta_t = \Phi_t - |\Omega_0|^{-1}\, \int_{\Omega_0} \Phi_t \,,
\end{equation}
so that $\int_{\Omega_0} \Theta_t = 0$, and $d\Theta_t = d\Phi_t$. Furthermore, the function $\Theta_t$ is $D$-symmetric.\medskip

Then,
\begin{equation}\label{E-jnn-42b}
\int_{\Omega_0} \Theta_t^2 = \int_{\Omega_0} \Phi_t^2 - |\Omega_0|^{-1} \, \left( \int_{\Omega_0} \Phi_t \right)^2\,.
\end{equation}

We introduce the notation,
$$
\delta(t) = d_r(\Omega_t,\Omega_0)\,.
$$

The constants $C_{25,i}$ which appear below only depend on $M$ and $s_0$.\medskip

Since $\int_{\Omega_t} \phi_t = 0$, we conclude as in the proof of Lemma~\ref{L-jnn-22s} that there exist constants $C_{25,1}$ and $C_{25,2}$ such that,
\begin{equation}\label{E-jnn-42c}\left\{
\begin{array}{l}
\big| \int_{\Omega_0} \Phi_t\big| \le C_{25,1}\, \delta(t) \,,\\[5pt]
\big| \int_{\Omega_0} \Phi_t^2 - 1 \big| \le C_{25,2}\, \delta(t) \,.
\end{array}
\right.
\end{equation}

 Using the condition \eqref{E-jnn-3}  to control $|\Omega_0|$, it follows that there exist constants $C_{25,3}$,\ldots,  $C_{25,5}$, such that
\begin{equation}\label{E-jnn-42d}
1 - C_{25,3}\, \delta(t)  \le \int_{\Omega_0} \Theta_t^2  \le 1 + C_{25,3}\, \delta(t)\,,
\end{equation}
and,  using Lemma~\ref{L-jnn-22s},
\begin{equation}\label{E-jnn-42e}
\int_{\Omega_0} |d\Theta_t|^2 \le \lambda_t + C_{25,4}\, \delta(t)^{q_0}
\le  \lambda_0 + C_{25,5}\, \delta(t)^{q_0}\,.
\end{equation}

Define the function
\begin{equation}\label{E-jnn-42f}
\Sigma_t := \Theta_t - \big( \int_{\Omega_0} \Theta_t \phi_0 \big) \Phi_0\,.
\end{equation}

Then $\Sigma_t$ is $D$-symmetric and satisfies
\begin{equation}\label{E-jnn-42g}
\int_{\Omega_0} \Sigma_t = 0  \text{~and~} \int_{\Omega_0} \Sigma_t \phi_0 = 0\,.
\end{equation}

It follows from our assumptions and notation that,
\begin{equation}\label{E-jnn-42h}
\int_{\Omega_0} |d\Sigma_t|^2 \ge \mu_0 \int_{\Omega_0} \Sigma_t^2\,,
\end{equation}
\begin{equation}\label{E-jnn-42k}
\int_{\Omega_0} \Sigma_t^2 = \int_{\Omega_0} \Theta_t^2 - \big( \int_{\Omega_0}\Theta_t \phi_0  \big)^2\,.
\end{equation}
 Using the fact that $(d\Phi_0)|_{\Omega_0} = d\phi_0$, and the variational definition of $(\lambda_0,\phi_0)$, we also have
\begin{equation}\label{E-jnn-42m}
\int_{\Omega_0} |d\Sigma_t|^2 = \int_{\Omega_0} |d\Theta_t|^2 - \lambda_0 \big( \int_{\Omega_0}\Theta_t \phi_0  \big)^2 \,.
\end{equation}

 From \eqref{E-jnn-42m} and the estimates on $\Theta_t$, there exists a constant $C_{25,6}$ such that
\begin{equation}\label{E-jnn-42n}
\int_{\Omega_0} |d\Sigma_t|^2 \le \lambda_0 - \lambda_0 \big( \int_{\Omega_0}\Theta_t \phi_0  \big)^2 + C_{25,6} \, \delta(t)\,.
\end{equation}

From \eqref{E-jnn-42h}, \eqref{E-jnn-42k} and  \eqref{E-jnn-42n}, it follows that there exist constants such that
\begin{equation}\label{E-jnn-42p}
\int_{\Omega_0} |d\Sigma_t|^2 \ge \mu_0 \left\{ 1 - C_{25,7} \, \delta(t) -   \big( \int_{\Omega_0}\Theta_t \phi_0  \big)^2 \right\}\,.
\end{equation}
and hence
\begin{equation}\label{E-jnn-42q}
\big| \int_{\Omega_0}\Theta_t \phi_0  \big| \ge 1 - \frac{(1+\mu_0)\,C_{25,8}}{\mu_0-\lambda_0}\, \delta(t)^{q_0}\,.
\end{equation}

 From \eqref{E-jnn-42q}, we deduce that for $\delta(t)$ small enough, the integral $\int_{\Omega_0}\Theta_t \phi_0$ is not zero.  Note that $\int_{\Omega_0} \Theta_t \phi_0 = \int_{\Omega_0} \Phi_t \phi_0$.  Write
$$
\int_{\Omega_0}\Phi_t \phi_0  = \int_{\Omega_0\cap \Omega_t}\phi_t \phi_0  + \int_{\Omega_0\setminus \Omega_t}\Phi_t \phi_0\,,
$$
and note that the second term tends to zero with $\delta(t)$. It follows that $\int_{\Omega_0 \cap \Omega_t} \phi_t \phi_0 \not = 0\,$, provided that $\delta(t)$ is small enough. This means that we can choose the sign of $\phi_t$ such that $\int_{\Omega_0\cap \Omega_t}\phi_t \phi_0 > 0\,$, provided that $\delta(t)$ is small enough. This proves the first assertion.\medskip

\emph{Proof of Assertion (2).}~ We now assume $\delta(t)$ to be small enough, so that we can uniquely determine the eigenfunction $\phi_t$ by $\|\phi_t\|_{L_2(\Omega_t)}=1$, with  $\int_{\Omega_0 \cap \Omega_t} \phi_t \phi_0 > 0$. More precisely, by \eqref{E-jnn-42q}, there exists a constant $C_{25,9}(M,s_0,\lambda_0,\mu_0)$ such that
\begin{equation}\label{E-jnn-44a}
\int_{\Omega_0}\Phi_t \phi_0 \ge 1 - C_{25,9}\, \delta(t)^{q_0}\,.
\end{equation}

Using \eqref{E-jnn-42c}, \eqref{E-jnn-44a}, and the fact that $\phi_0$ is normalized, there exists a constant $C_{25,10}(M,s_0,\lambda_0,\mu_0)$ such that
\begin{equation}\label{E-jnn-44c}
\int_{\Omega_0}(\Phi_t - \phi_0)^2 \le C_{25,10}\, \delta(t)^{q_0}\,.
\end{equation}

It follows that the functions $\Phi_t$ tend to $\phi_0$ in $L_2(\Omega_0)$.\medskip

 The family $\{\phi_t, t \ge 0\}$ is uniformly bounded in the $H^2(\Omega_t)$ (Proposition~\ref{P-jnn-6}), and hence the family $\{\Phi_t, t \ge 0\}$ is uniformly bounded in $H^2(\R^2)$, with compact support in $B(2M)$  (Proposition~\ref{P-jnn-8}). It follows that it is relatively compact in $C^{0,s_0-1}(\R^2)$, and weakly compact in $H^2(\R^2)$. The second assertion follows.\medskip

\emph{Proof of Assertion~(3).}~ Let $k$ be an integer, and let $K \subset \Omega_0$ be any compact subset. For $t$ small enough, we have $K \subset \Omega_t$. By interior regularity, $\Phi_t|_K = \phi_t|_K$ is uniformly bounded in $C^{k+1}(K)$ norm, and hence admits a convergent subsequence $\Phi_{t_j}$ in $C^k(K)$. Inequality \eqref{E-jnn-44c} shows that the limit of this subsequence must be $\phi_0$.  It follows that $\Phi_{t_j}$ converges to $\phi_0$ in $C^k(K)$. Because the limit is independent of the subsequence, it follows that $\phi_t$ tends to $\phi_0$ in $C^k(K)$. \hfill \qed \medskip

\textbf{Remark.~} Here is an alternative argument for the last assertion, which gives a stronger control of the convergence.\medskip

Let $\chi_1 ,\chi_2\in C_0^\infty (\Omega_0)$ such that $\chi_2=1$ on  $\mathrm{supp}(\chi_1) $.
We have
\begin{equation}\label{identite}
\begin{array}{ll}
 \Delta \chi_1 (\phi_t-\phi_0) & = [\Delta, \chi_1] (\chi_2 (\phi_t-\phi_0))  -  \chi_1 (\lambda_t \phi_t -\lambda_0 \phi_0)\\[5pt]
 &  = [\Delta, \chi_1] (\chi_2 (\phi_t-\phi_0))  - \chi_1 \lambda_t (\phi_t-\phi_0)\\[3pt]
& ~~~- \chi_1 (\lambda_t-\lambda_0) \phi_0\,.
 \end{array}
\end{equation}
Applying $(I-\Delta)^{-\frac 12}$ to this equality, and using Lemma \ref{L-jnn-22s} and \eqref{E-jnn-44c}, we get
$$
||\chi_1 (\phi_t-\phi_0)||_{H^1} \leq C \delta(t)^{\frac{q_0}{2}}\,.
$$
Hence, for any compact $K\subset \Omega_0$, we have
$$
|| \phi_t-\phi_0||_{H^1 (K)} \leq C \delta(t)^{\frac{q_0}{2}}\,.
$$
Similarly, starting from \eqref{identite}, given any $k \in \N$, and any compact $K$, we obtain,
$$
|| \phi_t-\phi_0||_{H^k (K)} \leq C(k,K) \delta(t)^{\frac{q_0}{2}}\,.
$$

\section{Domains with the symmetry of an equilateral triangle}\label{S-teqal}

\subsection{Preparation}\label{SS-teqal-1}

Let $\cT_e$ be the equilateral triangle, with vertices at $(0,0)$, $(1,0)$ and $(\frac{1}{2},\frac{\sqrt{3}}{2})$. The symmetry group of $\cT_e$ is generated by the mirror symmetries with respect to the side bisectors.
%% See Table~\ref{T-teqa-g0}.
\medskip

Up to scaling, the positive first Dirichlet eigenfunction of $\cT_e$ is given by the formula (see \cite{BH-ecp1}),
\begin{equation}\label{E-teqa-4d}\xi_1^{\mf{d}}(x,y):=
\sin(\frac{4\pi y}{\sqrt{3}}) + \sin\left( 2\pi (x - \frac{y}{\sqrt{3}}) \right) - \sin\left( 2\pi (x + \frac{y}{\sqrt{3}}) \right),
\end{equation}
which can also be written
\begin{equation}\label{E-teqa-4d2}
  \xi_1^{\mf{d}}(x,y) = 4 \sin\left( \frac{2\pi y}{\sqrt{3}} \right) \,
\sin\left( \pi (x - \frac{y}{\sqrt{3}}) \right) \, \sin\left( \pi (x + \frac{y}{\sqrt{3}}) \right).
\end{equation}

\begin{proposition}\label{P-teqa-0d}
The function $\xi_1^{\mf{d}}$ is positive in the interior of $\cT_e$. It has a unique critical point at $(\frac{1}{2},\frac{\sqrt{3}}{6})$, the centroid of the triangle. For $0 < c < \max_{\cT_e} \xi_1^{\mf{d}}$, the level curves $\{\xi_1^{\mf{d}} = c\}$ are smooth  strictly convex curves which have the same symmetries as $\cT_e$.
\end{proposition}%

This proposition is a consequence of \cite[Corollary~4.6]{CaFr1985}. We give an elementary proof using the following lemma \cite{Kaw1984}.

\begin{lemma}\label{L-teqa-0}
Let $\Omega$ be a convex bounded open set in $\R^2$. Let $\varphi$ be a positive, superharmonic function ($\Delta \varphi < 0$) in $\Omega$. If $\det \mathrm{Hess}\left( \log(\varphi) \right)$, the determinant of the Hessian of the function $\log (\varphi)$, is positive, then the super-level sets $\{\varphi > c\}$ are (strictly) convex.
\end{lemma}%

\emph{Proof of the lemma.} Let $\eta := \log(\varphi)$. Then,
$$
\varphi^2 \, \Delta \eta = \varphi \, \Delta \varphi - |d\varphi|^2\,.
$$
Since $\varphi$ is positive and superharmonic, it follows that $\Delta \eta < 0$, so that $\mathrm{Hess}(\eta)$ has at least one negative eigenvalue. On the other hand, since we work in dimension $2$, the positivity of $\det \mathrm{Hess}(\eta)$ implies that both eigenvalues of $\mathrm{Hess}(\eta)$ have the same sign. It follows that both eigenvalues are negative, and hence that $\mathrm{Hess}(\eta)$  is negative definite. The function $\varphi$ is (strictly) log-concave, and the lemma follows. \hfill \qed \medskip

\emph{Proof of the proposition.} It is easy to see that the only critical points of the function $\xi_1^{\mf{d}}$ in the closed triangle are the vertices and the centroid. This function is invariant under the  mirror symmetries with respect to the side bisectors of the triangle, and under the rotations with center the centroid, and angles $\pm \frac{2\pi}{3}$. It follows that its level sets have the same symmetries. Clearly, $\xi_1^{\mf{d}}$ is positive and superharmonic. It remains to show that $\det \mathrm{Hess}(\log \xi_1^{\mf{d}})$ is positive. This can be done by brute force. Let $\xi := \log (\xi_1^{\mf{d}})$. A Maple-aided computation gives,
\begin{equation}\label{E-teqa-P0d}
\left\{
\begin{array}{l}
\det \mathrm{Hess}(\xi) = \frac{4\pi^4}{3} \, \frac{N(\xi)}{D(\xi)}\,, \text{~with}\\[8pt]
N(\xi) = 2 - 2 \, \cos(\frac{2\pi y}{\sqrt{3}}) \, \cos\left( \pi (x - \frac{y}{\sqrt{3}})\right) \, \cos\left( \pi (x + \frac{y}{\sqrt{3}})\right)\,,\\[8pt]
D(\xi) = (\xi_1^{\mf{d}})^2\,.
\end{array}
\right.
\end{equation}

The proof of Proposition~\ref{P-teqa-0d} is complete.\hfill \qed \bigskip

\textbf{Notation}. We shall now work with the equilateral triangle $\cT_0$, with vertices $A = (-\frac{1}{2}, -\frac{\sqrt{3}}{6})$, $B = (\frac{1}{2}, -\frac{\sqrt{3}}{6})$ and $C = (0,\frac{\sqrt{3}}{3})$, and centroid $O = (0,0)$. Making the change of coordinates $x=\frac 12 + u$ and $y = \frac{\sqrt{3}}{6} + v$, in $\xi_1^{\mf{d}}$, we obtain a first Dirichlet eigenfunction for $\cT_0$,
\begin{equation}\label{E-teqa-4d2-0}\small %\scriptstyle
\varphi_1^{\mf{d}}(u,v) = 4 \sin \frac{\pi}{3}(1+2\sqrt{3}v)\,
\sin \frac{\pi}{3}(1-3u+\sqrt{3}v)\,
\sin \frac{\pi}{3}(1-3u-\sqrt{3}v)\,.
\end{equation}

Define the function,
\begin{equation}\label{E-teqa-f0}
f_0(u,v) := (1+2\sqrt{3}v)\, (1+3u-\sqrt{3}v) \, (1-3u-\sqrt{3}v) \,.
\end{equation}

\begin{proposition}\label{P-teqa-0f}
The function $f_0$ is positive in the interior of $\cT_0$. It has a unique critical point at $O$, the centroid of the triangle. For $0 < c < 1$, the level curves $\{f_0 = c\}$ are smooth  strictly convex curves which have the same symmetries as $ \cT_0$.
\end{proposition}%

\proof We again make use of Lemma~\ref{L-teqa-0}. The first two assertions are clear. The function $f_0$ is clearly invariant under the symmetries of $\cT_0$, so are its level sets. An easy computation gives $\Delta f_0 = - 36$, so that $f_0$ is superharmonic. Let $g := \log(f_0)$.  Define the functions $A_{uu}, A_{uv}$ and $A_{vv}$ by the formulas
$$
A_{uu} = f_0^2 \, \frac{\partial^2 g}{\partial u^2}\,, \text{~etc.} \,.
$$
Then,
$$
f_0^4 \, \det \mathrm{Hess}(g) = A_{uu}\, A_{vv} - (A_{uv})^2\,.
$$

A Maple-aided computation gives,
$$
A_{uu}\, A_{vv} - (A_{uv})^2 = 324 \, f_0^2 \, (1+6u^2+6v^2)\,,
$$
so that
\begin{equation}\label{E-teqa-P0da}
\det \mathrm{Hess}(g)(u,v) =324 \, \frac{1+6u^2+6v^2}{f_0^2(u,v)}\,.
\end{equation}

This completes the proof of Proposition~\ref{P-teqa-0f}. \hfill \qed %\medskip

\begin{figure}
  \centering
  \includegraphics[scale=0.40]{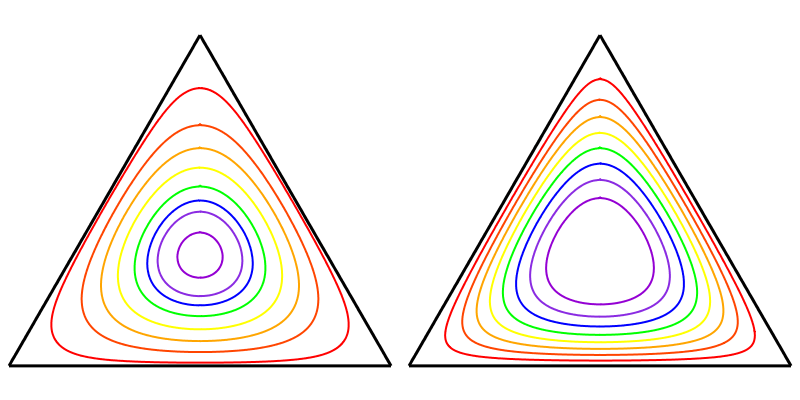}%
\caption{Level sets of $\varphi_1^{\mf{d}}$ (left) and $f_0$ (right)}\label{F-TLevels}
\end{figure}

\FloatBarrier

\begin{remark}\label{R-teqa-0dg}
Note that the function $f_0$ is (up to scaling) the torsion (or warping) function of the equilateral triangle, see \cite[Section~7]{HelPer2016}. The square root of the warping function $f_{\Omega}$ is known to be strictly concave, see \cite[Theorem~4.1]{Ken1985}
 \end{remark}%

\subsection{Domains with $\cG_0$-symmetry}\label{SS-teqal-2}

\begin{figure}[!ht]
  \centering
  \includegraphics[scale=0.40]{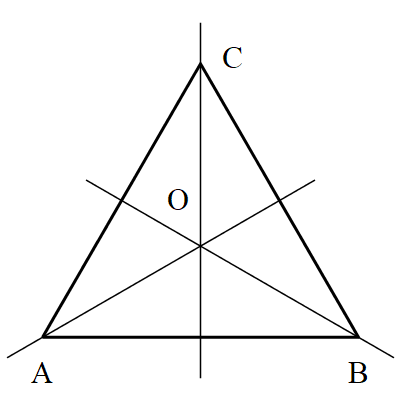}
\caption{The equilateral triangle $\cT_0$}\label{F-teq0}
\end{figure}

Recall that $\cT_0$ is the equilateral triangle with vertices $A$, $B$, and $C$, and centroid $O$. Call $D_A$, $D_B$ and $D_C$ the bisectors of its sides. The coordinates are chosen so that $
 D_C = \{u=0\}$, see Figure~\ref{F-teq0}.\medskip

The isometry group of $\cT_0$ is the group
\begin{equation}\label{E-teqa-2}
\cG_0 = \left\lbrace  I, D_A, D_B, D_C, R, R^2 \right\rbrace \,.
\end{equation}
where $D_A$, is the mirror symmetry with respect to the bisector $D_A$, $R$ the rotation with center $0$ and angle $\frac{2\pi}{3}$.
% See Table~\ref{T-teqa-g0}.}
%\medskip

\begin{figure}
  \begin{center}
  \includegraphics[scale=0.8]{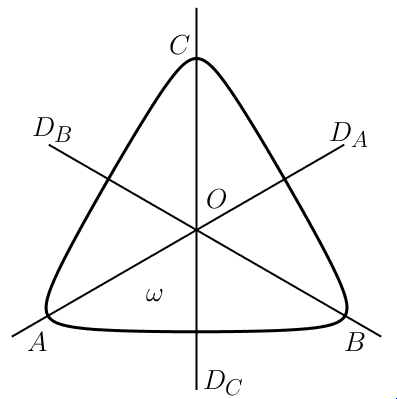}
  %vspace{-105mm}
\caption{A domain $\Omega$ in the class \eqref{E-teqa-4}}\label{F-teqa-sym}
\end{center}
\end{figure}

To construct smooth counterexamples to $\ecp$, the idea is to start from the equilateral triangle, and to consider  the class $ \cL_{M,0}$ of domains $\Omega$ with the following properties,
\begin{equation}\label{E-teqa-4}
\left\{
\begin{array}{l}
\Omega \in \cL_M\,,\\[5pt]
%\Omega \text{~has Lipschitz constant} \le M\, \\[5pt]
\Omega \text{~admits~} \cG_0 \text{~as symmetry group,}
\end{array}%
\right.
\end{equation}
see Figure~\ref{F-teqa-sym}.\medskip

The super-level sets $\left\lbrace x \in \cT_0 ~|~  \varphi_1^{\mf{d}} > c \right\rbrace$ of the first Dirichlet eigenfunction, and the super-level sets $\left\lbrace x \in \cT_0 ~|~  f_0 > c \right\rbrace$ of the torsion function $f_0$ provide examples of $C^{\infty}$  strictly convex domains $\Omega$ with the symmetry group $\cG_0$, see Figure~\ref{F-TLevels}. Another example is the equilateral triangle with rounded corners, $\cT_{0,a}$: replace each corner by an arc of circle, with radius $a$, centered on the corresponding bisector, and tangent to the sides. This yields a convex domain, with $C^1$, piecewise $C^2$, boundary, with symmetry group $\cG_0$.\medskip

One can show that these families of domains belong to the class $\cL_{M,0}$ for some $M > 0$, see \eqref{E-teqa-4}. In Section~\ref{S-appl},  in order to prove Theorem~\ref{T-intro-4}, we shall consider yet another family, and prove that it is indeed in the class $\cL_M$ for some $M$.\medskip

We conclude this section with a spectral property of the domains in the class $\cL_{M,0}$.

\begin{proposition}\label{P-teqa-n2}
Let $\Omega$ be a smooth domain in the class \eqref{E-teqa-4}. Then, the first Neumann eigenvalues of $\Omega$ satisfy
\begin{equation}\label{E-teqa-n6}
0 = \nu_1 < \nu_2 = \nu_3 < \nu_4 \le \cdots
\end{equation}
More precisely, the eigenspace $\cE(\nu_2)$ admits a basis $\{\phi,\psi\}$ such that $D_C^{*}\phi = \phi$, and $D_C^{*}\psi = -\psi$. Furthermore, $\cZ(\phi) \cap D_C = \{O\}$, and $ \cZ(\psi) = D_C \cap \Omega$.
\end{proposition}%

\proof
The proof is based on the following properties: a Neumann eigenfunction $\xi$ of $\Omega$ has finitely many interior critical zeros, finitely many boundary zeros, and its nodal set consists of finitely many simple regular arcs whose end points are either interior critical zeros, or boundary zeros. We do not need to know the local structure at boundary zeros.\medskip

Let $\xi \not = 0$ be a 2nd Neumann eigenfunction. Assertions (a)--(c) hold for a simply-connected regular domain.  \medskip

\textbf{(a)} \emph{The nodal set $\cZ(\xi)$ cannot contain any interior closed curve.} Indeed, there would otherwise exist an interior nodal domain $\omega$ of $\xi$, for which we could write $\nu_2(\Omega) = \delta_1(\omega) > \delta_1 (\Omega)$, contradicting the inequality $\nu_2(\Omega) < \delta_1(\Omega)$ due to to P\'{o}lya \cite{Pol1952} and Szeg\"{o} \cite{Sze1954} (here the $\delta$'s refer to Dirichlet eigenvalues).\medskip

\textbf{(b)} \emph{The nodal set $\cZ(\xi)$ does not contain any interior critical zero.} Assume this is not the case. Then, there would exist an interior critical $x_0$, and at least four semi-arcs issuing from $x_0$ and contained in $\cZ(\xi)$. Continuing these arcs, we either obtain a closed nodal curve, or reach the boundary at distinct points. The first case is impossible by (a). In the second case, because $\Omega$ is simply-connected, we would obtain at least four nodal domains, contradicting Courant's theorem.\medskip

\textbf{(c)} \emph{The nodal set of any 2nd Neumann eigenfunction $\xi$ in $\Omega$ consists of a single simple curve without critical zeros, meeting the boundary at two distinct points.} The fact that such a curve must be simple and without critical zeros follows from (a) and (b). The fact that its boundary points must be distinct follows from (a). Assume that there exist two such curves. By (b), they cannot meet in the interior of $\Omega$. If they had identical boundary points, we would get a contradiction by (a). In the other case, we would get a contradiction with Courant's theorem.\medskip

We now assume that $\Omega$ has the symmetries of the equilateral triangle.\medskip

Let $D := D_C$, and define the set of functions
\begin{equation}\label{E-teqa-n10}
\cS_{\sigma} := \left\lbrace  \varphi ~|~ D^* \varphi = \sigma \varphi \right\rbrace \,, ~~ \sigma \in \{+,-\}\,.
\end{equation}

Because $D$ is an isometry, $D^*$ leaves $\cE(\nu_2)$ globally invariant, and the eigenspace decomposes as
\begin{equation}\label{E-teqa-n12}
\cE(\nu_2) = ( \cE(\nu_2) \cap \cS_{+})  \oplus (\cE(\nu_2) \cap \cS_{-}) \,.
\end{equation}

Because the rotation $R$ is an isometry, $R^*$ leaves $\cE(\nu_2)$ globally invariant, and so does the map
\begin{equation}\label{E-teqa-n14}
T := R^* - R^{*2}
\end{equation}
which commutes with the Laplacian $\Delta$. \medskip

It is easy to see that $D^* \circ T = - T \circ D^*$, so that
\begin{equation}\label{E-teqa-n16}
T \left(  \cE(\nu_2) \cap \cS_{\pm} \right) \subset \cE(\nu_2) \cap \cS_{\mp} \,,
\end{equation}
\begin{equation}\label{E-teqa-n18}
\ker (T) = \ker (R^*-I)\,,
\end{equation}
and that
\begin{equation}\label{E-teqa-n19}
\cS_{\sigma} \cap \ker (T) = \big\lbrace  \varphi ~|~ D_M^{*}\varphi = \sigma \, \varphi\,, ~\forall M \in \{A,B,C\} \big\rbrace \,.
\end{equation}

The following assertions hold.\medskip

\textbf{(d) }\emph{If $0 \not = \xi \in \cE(\nu_2)$, then $R^{*}\xi \not = \xi$}. Indeed, using (c) and the $R$-invariance of $\xi$, $\cZ(\xi)$ would contain at least three boundary points, contradicting (c).\medskip

\textbf{(e)} \emph{The dimension of $\cE(\nu_2)$ is at least $2$.} Indeed, we would otherwise have $\dim \cE(\nu_2) = 1$, and hence, for some $0 \not = \xi \in \cE(\nu_2)$, $R^{*}\xi = \pm \xi$. Since $R^{*}=I$, this would imply that $R^{*}\xi=\xi$, contradicting (d).\medskip

\textbf{(f)} \emph{The dimension of $\cE(\nu_2) \cap \cS_{-}$ is at most $1$.} Indeed, if $0\not = \xi \in \cE(\nu_2) \cap \cS_{-}$, then $\xi$ vanishes on $D\cap \Omega$, and it cannot vanish elsewhere by Courant's theorem. This implies that $\xi|_{\Omega_+}$ is the first eigenfunction of $\Omega_{+}$ (with mixed boundary conditions), and hence unique up to scaling. This implies that $\xi$ itself is unique up to scaling.\medskip

\textbf{(g)} \emph{The dimension of $\cE(\nu_2) \cap \cS_{+}$ is at least $1$.} Indeed, by (e) and (f), there exists $0 \not = \xi \not \in \cE(\nu_2) \cap \cS_{-}$. This implies that $\phi := \frac{1}{2}(\xi + D^{*}\xi)$ is a nonzero function in $\cE(\nu_2) \cap \cS_{+}$.\medskip

\textbf{(h)} \emph{Both spaces $\cE(\nu_2) \cap \cS_{\pm}$ have dimension $1$, and there exists a basis $\{\phi,\psi\}$ of $\cE(\nu_2)$, such that $\phi$ is $D$-symmetric, and $\psi$ $D$-anti-symmetric.} Using \eqref{E-teqa-n14}, we see that $T(\xi)=0$ if and only if $R^{*}=\xi$, so that $T$ is injective from $\cE(\nu_2)$ into itself. Using \eqref{E-teqa-n16} in both directions, we infer that
$\dim \cE(\nu_2) \cap \cS_{\pm} = 1$, and the assertion follows.\medskip

\textbf{(i)} \emph{We have $\cZ(\psi) = D\cap \Omega$ and $\cZ(\phi) \cap D = \{O\}$.} We have already proved the first part of the assertion in (f). Up to scaling, we have $\phi = T(\psi)$. Since $R(O)=O$, the definition of $T$ implies that $\phi(O)=0$. The fact that $\cZ(\psi)$ meets $D$ at exactly one point follows from (a), (c) and the symmetry of $\phi$.\hfill \qed

\begin{remark}\label{R-teqa-n2R}
Note that the inequality $\nu_2(\omega) < \delta_1(\omega)$ is valid for any sufficiently regular, bounded domain, without any convexity assumption. The fact that a second Neumann eigenfunction cannot have a closed nodal line motivated the ``closed nodal line conjecture for a second Dirichlet eigenfunction'', see \cite{Pay1967}, last paragraph on page 466, and Conjecture~5, and \cite{Ken2017}.
\end{remark}%

\begin{remarks}\label{R-teqa-n6} Concerning the multiplicity of $\nu_2$, we can mention the following.
\begin{enumerate}
  \item According to \cite[Remarks~(2), p. 206]{LeWe1986}, if $\Omega$ is close enough to $\cT_0$ in the sense of the Hausdorff distance, then $$\dim \cE\left( \nu_2(\Omega)\right) = \dim \cE\left( \nu_2(\cT_0) \right) = 2\,.$$
  \item For any smooth simply-connected domain $\Omega$, $\dim \cE\left( \nu_2(\Omega)\right) \le 3$. This bound was first given by Cheng \cite{Che1976} for smooth simply-connected compact surfaces without boundary, see also the assertion in \cite[line (-8), p.~1170]{HoMiNa1999}. In this latter paper, the authors indicate that the assumption that $\Omega$ is smooth is probably too strong. The smoothness assumption is used to describe the local behaviour of the nodal set at a boundary point. In the non-smooth case, it might be possible to obtain a result on the local structure of the nodal set similar to the one described by Alessandrini \cite{Ale1994} for the Dirichlet boundary condition. See also \cite{HHOT2009}
  \item In \cite[Theorem~2.3]{Lin1987}, Lin proved that the second Dirichlet eigen\-space of a smooth convex domain has dimension at most $2$.
\end{enumerate}
\end{remarks}%
\medskip

\section{Proof of Theorem~\ref{T-intro-4}}\label{S-appl}
 To prove Theorem~\ref{T-intro-4}, we apply the deformation technique of Section~\ref{S-jnn} to a special family of domains in the class $\cL_{M,0}$.

%% \begin{equation}\label{E-jnn-46}
%% \cL_{M,0} := \left\{ \Omega \in \cL_M ~|~ \Omega \text{~admits~}
%% \cG_0 \text{~as symmetry group} \right\}\,.
%% \end{equation}

\subsection{Construction of the family $\Omega_t$}

Let $t$ be a nonnegative parameter. Introduce the function
\begin{equation}\label{E-appln-2}
f_{0,t}(u,v) = (1+t+2\sqrt{3}v)\, (1+t+3u-\sqrt{3}v) \, (1+t-3u-\sqrt{3}v) \,.
\end{equation}

When $t=0$, we recover the function $f_0$ defined by \eqref{E-teqa-f0}. When $t > 0$, the function $f_{0,t}$ is a torsion function for the triangle $\cT_{0,t}$ obtained from $\cT_0$ by dilation of ratio $(1+t)$. This equilateral triangle has vertices
$A_t = (-\frac{1+t}{2}, -\frac{\sqrt{3}(1+t)}{6}), ~~B = (\frac{1+t}{2}, -\frac{\sqrt{3}(1+t)}{6})$ and $C = (0,\frac{\sqrt{3}(1+t)}{3})$.\medskip

An immediate computation gives that
\begin{equation}\label{E-appln-2c}
f_{0,t}(A) = f_{0,t}(B) = f_{0,t}(C) = t^2 (3+t)\,.
\end{equation}

\begin{definition}\label{D-appln-2}
Define the domain $\Omega_t$ to be the super-level set
\begin{equation}\label{E-appln-2e}
\Omega_t := \left\lbrace f_{0,t} > t^2 (3+t) \right\rbrace\,.
\end{equation}
\end{definition}%

\begin{figure}[!ht]
  \centering
  \includegraphics[scale=0.40]{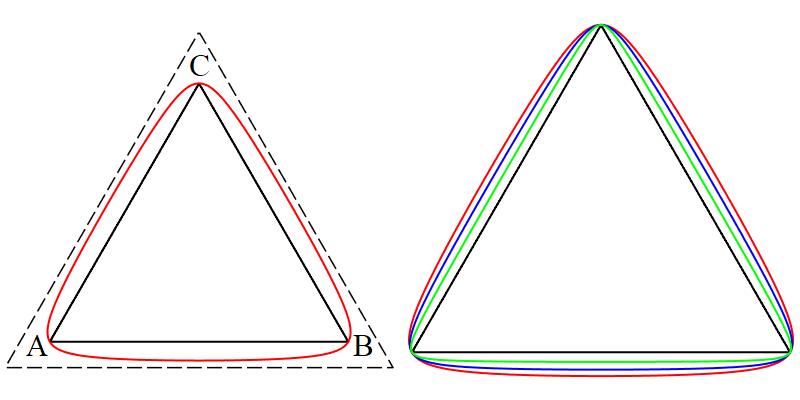}
  \vspace{-5mm}
\caption{Domains $\Omega_t$}\label{F-appln-2}
\end{figure}

The triangle $\cT_0$, the triangle $\cT_{0,t}$ (dashed line), and a domain $\Omega_t$ (red line) are displayed in Figure~\ref{F-appln-2}, left. The triangle $\cT_0$, and domains $\Omega_t$, with $t=0.3$ (red), $t=0.2$ (blue), and $t=0.1$ (green), are displayed in Figure~\ref{F-appln-2}, right.\medskip

\FloatBarrier

Let us summarize the properties of the domains $\Omega_t$.

\begin{proposition}\label{P-appln-2}
The family of domains $\{\Omega_t\}_{0\le t \le \frac{1}{2}}$ has the following properties.
\begin{enumerate}
  \item $\Omega_0 = \cT_0$.
  \item  For $t > 0$, the domain $\Omega_t$ is strictly convex, bounded, and open, with $C^{\infty}$ boundary. Furthermore, $\cT_0 \subset \Omega_t$, and $A, B, C \in \partial \Omega_t$.
  \item The domain $\Omega_t$ has the symmetry group $\cG_0$.
  \item The family $\Omega_t$ is increasing, for $0 < t_1 < t_2$,
  $$\Omega_{t_1} \subset \Omega_{t_2}\,.$$
  \item For $0 \le t \le \frac{1}{2}$, the domains $\Omega_t$ belong to the class $\cL_M$ for some positive constant $M$.
\end{enumerate}
\end{proposition}%

\begin{figure}[!ht]
  \centering
  \includegraphics[scale=0.50]{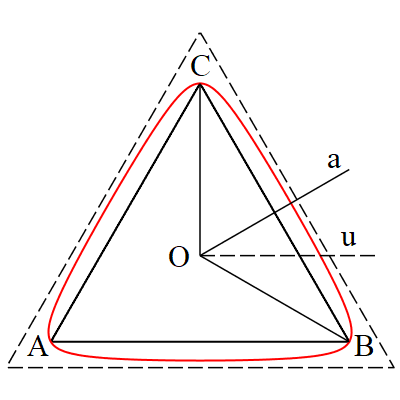}
  \vspace{-5mm}
\caption{}\label{F-appln-4}
\end{figure}

\proof Assertion~(1) is obvious.\\[3pt]
Assertion~(2).  The first part follows from Proposition~\ref{P-teqa-0f} by dilation of ratio $(1+t)$.  For the second part, note that, by definition of $\Omega_t$, the vertices $A, B$ and $C$ belong to $\partial \Omega_t$. The inclusion (of open sets) $\cT_0 \subset \Omega_t$ follows from the convexity of $\Omega_t$.\\[3pt]
Assertion~(3). This follows from Proposition~\ref{P-teqa-0f}.\\[3pt]
Assertion~(4). The domain $\Omega_t$ can also be defined by $\{g_t > 0\}$, where %
$$
g_t(u,v) = f_{0,t}(u,v) - t^2(3+t) = f_0(u,v) + 3t - 9t(u^2+v^2)\,.
$$
Let $t_1 < t_2$. To prove that $\Omega_{t_1} \subset \Omega_{t_2}$, it suffices to consider the points $(u,v) \in \Omega_{t_1} \sm \cT_0$. For such $(u,v)$, we have $g_{t_1}(u,v) > 0$ and $f_0(u,v) \le 0$. This implies that
$$
3t_1 (1-3u^2-3v^2) > -f_0(u,v) \ge 0\,,
$$
and hence that $1-3u^2-3v^2 > 0$. On the other hand, we have
$$
g_{t_1}(u,v) - g_{t_2}(u,v) = 3(t_1-t_2)(1-3u^2-3v^2) < 0\,,
$$
i.e., $g_{t_2}(u,v) > 0$, or $(u,v) \in \Omega_{t_2}$.\\[3pt]
Assertion~(5). Since $\cT_0 \subset \Omega_t \subset \cT_{0,t}$, the domains satisfy condition~\eqref{E-jnn-3}. It remains to show that they satisfy condition~\eqref{E-jnn-6}, i.e., that they can be defined in polar coordinates, as
$$
\Omega_t = \{(r,\theta) ~|~ 0 \le r < \rho(t,\theta)\}
$$
where the functions $\rho(t,\cdot)$ are uniformly Lipschitz. Due to rotational invariance, it suffices to look at the part of $\partial \Omega_t$ contained in the sector $BOC$, see Figure~\ref{F-appln-4}. This part of the boundary is symmetric with respect to the bisector $D_A$, so that it suffices to look at the sector $BOa$. With respect to the $u$-axis $Ou$, the angle $\theta$ then varies from $-\frac{\pi}{6}$ ($OB$) to $\frac{\pi}{6}$ ($Oa$). \medskip

Instead of polar coordinates $(\rho,\theta)$, we use ``inverse'' polar coordinates $(r,\theta)$, where $r \, \rho \equiv 1$. The inverse polar equation of the side $BC$ of $\cT_0$, is
\begin{equation}\label{E-appln-4a}
r_A(\theta) = 2 \sqrt{3}\cos(\theta-\frac{\pi}{6})\,, \text{~~~for~} \theta \in [-\frac{\pi}{6},\frac{\pi}{6}]\,.
\end{equation}

Let $r=r(t,\theta)$ be the inverse polar equation of the arc $BC \subset \partial \Omega_t$. Because $\cT_0 \subset \cT_{0,t}$, we have
\begin{equation}\label{E-appln-4c}
\frac{1}{1+t}\, r_A(\theta) \le r(t,\theta) \le r_A(\theta) \text{~for~}
\theta \in [-\frac{\pi}{6},\frac{\pi}{6}]\,.
\end{equation}

Using the definition of $\Omega_t$, we also have that $r(t,\theta)$ is a root of the equation
\begin{equation}\label{E-appln-4e}
(1+3t)\, r^3 - 9(1+t)\, r + 6\sqrt{3}\sin^3(\theta) - 18 \sqrt{3} \sin(\theta)\cos^2(\theta)=0\,.
\end{equation}
or, equivalently,
\begin{equation}\label{E-appln-4f}
(1+3t)\, r^3 - 9(1+t)\, r - 6\sqrt{3}\sin(3\theta)=0\,.
\end{equation}

Looking at the global picture of $f_{0,t}^{-1}(0)$, it is easy to see that this equation has one simple root satisfying \eqref{E-appln-4c}. Taking the derivative $r_{\theta}$ with respect to $\theta$, we obtain,
\begin{equation}\label{E-appln-4h}
\big( (1+3t)\,r^2 - 3(1+t)\big) r_{\theta} -6\sqrt{3}\cos(3\theta) = 0\,.
\end{equation}

Note that
\begin{equation}\label{E-appln-4j}
(1+3t)\, r^3 - 3(1+t)\, r = \big( (1+3t)\,r^3 - 9(1+t)\, r \big) + 6(1+t)\, r\,,
\end{equation}
so that
\begin{equation}\label{E-appln-4k}
(1+3t)\, r^3 - 3(1+t)\, r = 6\big( (1+t)\,r + \sqrt{3}\sin(3\theta)\big)\,.
\end{equation}

Using \eqref{E-appln-4c}, we have
\begin{equation}\label{E-appln-4m}
(1+t)\,r + \sqrt{3}\sin(3\theta) \ge 2\sqrt{3}\cos(\theta-\frac{\pi}{6}) + \sqrt{3}\cos(3(\theta-\frac{\pi}{6}))\,,
\end{equation}
and hence
\begin{equation}\label{E-appln-4p}
(1+3t)\, r^3 - 3(1+t)\, r \ge 6 \sqrt{3}\cos(\theta-\frac{\pi}{6})\big( 4\cos^2(\theta-\frac{\pi}{6})-1\big)\,.
\end{equation}

It follows that $r_{\theta}$ is positive in the interval $]-\frac{\pi}{6}\,,\,\frac{\pi}{6}\,[$, and that
\begin{equation}\label{E-appln-4r}
0 \le \frac{r_{\theta}(t,\theta)}{r(t,\theta)} \le \tan(\frac{\pi}{6} -\theta) \le \sqrt{3}\,.
\end{equation}

Note that $r(t,\theta) \ge 2\sqrt{3}$. This proves that condition \eqref{E-jnn-6} is satisfied. \hfill \qed

\subsection{Proof of Theorem~\ref{T-intro-4}}\label{SS-appl3}

\begin{figure}[!ht]
  \centering
  \includegraphics[scale=0.60]{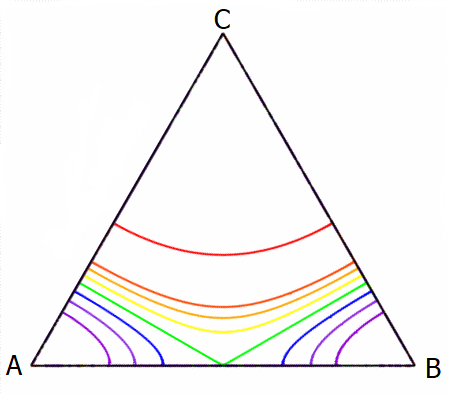}
\caption{Level lines of the second  symmetric Neumann eigenfunction of the equilateral triangle}\label{F-appl-6}
\end{figure}

The fact that the equilateral triangle $\cT_0$ provides a counterexample to $\ecp(\cT_0,\mf{n})$ follows from the analysis of the level lines of the $D$-symmetric second Neumann eigenfunction $\phi_{\cT_0}$, see  \cite[Section~3]{BH-ecp1}. Some of the levels lines of $\phi_{\cT_0}$ are displayed in Figure~\ref{F-appl-6}. \medskip

\begin{figure}[!ht]
  \centering
  \includegraphics[scale=0.70]{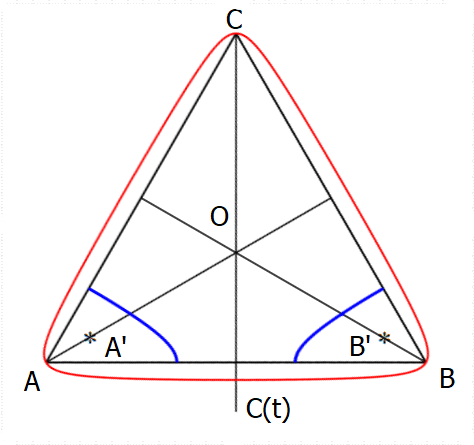}
\caption{Proof of Theorem~\ref{T-intro-4}.}\label{F-appl-8}
\end{figure}

Deform the domain $\cT_0 =: \Omega_0$ using Proposition~\ref{P-appln-2}. Denote the normalized $D$-symmetric eigenfunctions by $\phi_t$, and their extensions by $\Phi_t$.\medskip

The function $\phi_0$ is such that $\phi_0(C) > 0$, and $\phi_0(A)=\phi_0(B) < 0$, see Figure~\ref{F-appl-8}. According to \cite[Section~3]{BH-ecp1}, we now choose (and fix) some $a > 0$, such that $\{\phi_0 +a = 0\}$ consists of two disjoint arcs, symmetric with respect to the side bisector $D_C$ (blue arcs in the figure). We have $\phi_0 + a > 0$ in the connected component of $\cT_0\sm\{\phi_0 +a = 0\}$ which contains $O$, and $\phi_0 + a < 0$ in the two connected components close to the vertices $A$ and $B$. Choose $A'$ and $B'$ in these connected components. Note that $\phi_0 |_{D_C}  + a > 0$, and $\phi_0(A') + a =\phi_0(B') + a  < 0$.\medskip

We now consider the family $\Omega_t$. Apply Lemma~\ref{L-jnn-32s} to the family $\phi_t$, and get that for $t$ sufficiently small $$\phi_t(A')+ a =\phi_t(B') + a < 0\,.$$

Call $C(t)$ the intersection point of the bisector $D_C$ with $\partial \Omega_t$, opposite to the vertex $C$.\medskip

\textbf{Claim~1.~} For $t$ sufficiently small, $\phi_t|_{[CC(t)]} + a > 0$. \smallskip

Indeed, we could otherwise find a sequence $t_k$, tending to zero, and a point $m_k \in [CC(t_k)]$, such that $\phi_{t_k}(m_k)+a \le 0$.
The family $\Phi_{t_k}$ is bounded in $H^2$ with compact support in $B(0,2M)$.  Hence,  there exists a subsequence $t'_j$ which tends to $0$, and a function $\Phi \in C^{0}(\R^2)\cap H^2(\R^2)$ such that $m_k$ converges to  some $m \in [CC(0)]$ and $\Phi_{t'_j}$ converges to $\Phi$ uniformly in $B(2M)$, and in particular in $\wb{\cT_0}$. Since, by Lemma~\ref{L-jnn-32s},   $\Phi_{t'_j} |_{\cT_0}$ converges to $\phi_0$ in $\mathcal D'(\cT_0)$, it follows that $\phi_0 = \Phi|_{\cT_0}$ and this extends by continuity to $\wb{ \cT_0}$. In particular, we would get $\Phi(m) +a = \phi_0(m) + a \leq 0$. A contradiction.\medskip

The claim proves that  for $t$ small enough, the points $A'$ and $B'$ belong to distinct connected components of  $\Omega_{t }\sm\{\phi_{t}+a=  0\}$, so that $\phi_{t}+a$ has at least three connected component (a ``positive'' one, and two ``negative ones''). \medskip

In particular this proves that, for $t$ small enough, the domains $\Omega_t$ provide a counterexample to the Extended Courant property. \medskip

 We shall now prove that, for $t$ small enough, $\phi_t + a$ has exactly three nodal domains.

\begin{lemma}\label{L-Glad}
Let $\{\varphi_n, n\ge 1\}$ be an orthonomal basis of eigenfunctions of the Neumann problem in a bounded domain $\Omega$, associated with the eigenvalues $0 = \nu_1(\Omega) < \nu_2(\Omega) \le \ldots$. Choose $\varphi_1$ (a constant function) to be positive. Then, for any $a > 0$, the set $\Omega\sm\{\varphi_n+a\varphi_1 =0\}$ has at most $(n-1)$ connected components in which $\varphi_n+a$ is positive.
\end{lemma}%

\begin{remark}\label{R-appln-4}
A statement analogous to Lemma~\ref{L-Glad}, for the Dirichlet problem in $\Omega$, appears as Theorem~1 in \cite{GlZh2003}. The proof given by Gladwell-Zhu is similar to the proof of Courant's nodal domain theorem, and turns out to apply to both the Dirichlet and the Neumann boundary conditions, hence to Lemma~\ref{L-Glad}. The examples of rectangles with cracks in \cite[Section~3]{BH-ecp1} show that one can a priori not control the number of connected components of $\Omega\sm\{\varphi_n+a\varphi_1 =0\}$ in which $\varphi_n+a$ is negative.
\end{remark}%

We proceed with the proof that, for $t$ small enough, $\phi_t + a$ has exactly three nodal domains. According Lemma~\ref{L-Glad}, we have to prove that $\{\phi_t +a <0\}$ has at most two connected components. The proof goes as follows.\smallskip

\begin{figure}
  \centering
  \includegraphics[scale=0.4]{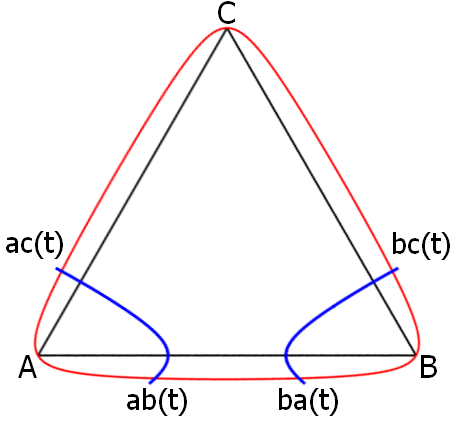}
  \caption{}\label{F-pthm-f2}
\end{figure}

First, we observe that $\phi_0$ is naturally defined as a trigonometric polynomial on all $\mathbb R^2$.  Observe that for $t$ small enough, $\{ \phi_0 + a =0 \} \cap \Omega_t$ consists of two symmetric curves crossing $\partial \Omega_t$ transversally at the points $ac(t), ab(t), ba(t),bc(t)$. As $t$ tends to $0$, these points tend to the intersection points of $\{\phi_0 +a=0\}$ with $\partial \cT_0$, see Figure~\ref{F-pthm-f2}.\medskip

 For $\epsilon >0$ small enough, we introduce,
\begin{equation}\label{E-appln-10}
\Omega_- (a+\epsilon,\phi_0,t):= \{\phi_0+a +\epsilon \leq  0\} \cap \Omega_t\,,
\end{equation}

\begin{equation}\label{E-appln-12}
\Omega_+ (a-\epsilon,\phi_0,t):= \{\phi_0+a-\epsilon \geq  0\} \cap \Omega_t\,,
\end{equation}
and
\begin{equation}\label{E-appln-13}
\Omega(a,\varepsilon,\phi_0) := \{-\varepsilon \le \phi_0+a \le \varepsilon\}\cap \Omega_t\,.
\end{equation}

\begin{figure}
  \centering
  \includegraphics[scale=0.4]{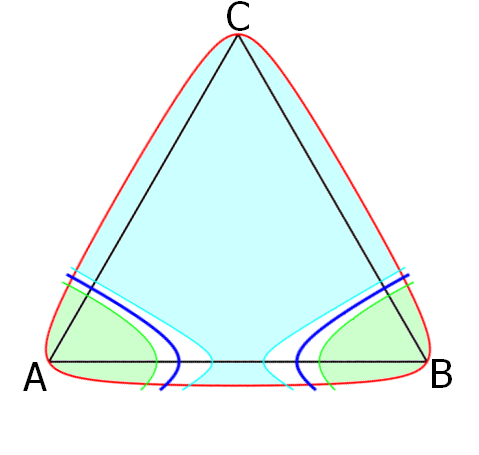}
  \caption{Localization}\label{F-pthm-f4}
\end{figure}

These domains are displayed respectively in green, blue, and white in Figure~\ref{F-pthm-f4}.\medskip

\textbf{Claim~2.~} For $t$ small enough,
\begin{equation}\label{E-appln-14}
\left\{
\begin{array}{ll}
\Omega_- (a+\epsilon,\phi_0,t) & \subset \{\phi_t+ a < 0\}\,,\\[5pt]
\Omega_+ (a-\epsilon,\phi_0,t) & \subset \{\phi_t+ a > 0\}\,.
\end{array}%
\right.
\end{equation}

Indeed, if the first inclusion were not true, there would exist a sequence $t_n>0$, tending to $0$, and $x_n\in \Omega_{t_n}$,  such that $\phi_{t_n} (x_n) + a \geq 0$ and $\Phi_{t_n}$ bounded in $H^2$. As above, after extraction of a subsequence we can assume that $x_n \rightarrow x_\infty$, and that $\Phi_{t_n}$ tends to $\Phi$ in $C^0$. This implies the existence of $x_\infty$ such that $\Phi (x_\infty)=\phi_0(x_\infty) \geq -a$. But $x_\infty \in  \Omega_- (a+\epsilon,\phi_0,0)$ leading to a contradiction. The second inclusion can be proved in a similar way.\medskip

As a consequence, for $t$ small enough, there are two symmetric components of $\{ \phi_t + a < 0\}$, each one containing a component of  $\{\phi_0+a+\epsilon \leq  0\} \cap \Omega_t$. Furthermore,  the ``positive'' component of $\phi_t +a$ contains  $\Omega_+ (a-\epsilon,\phi_0,t) $.\smallskip

 We deduce from this localization, that a third ``negative'' connected component of $\phi_t +a$,  if any, is necessarily contained in $\Omega(a,\varepsilon,\phi_0)$, hence stays away from the vertices $A$, $B$ and $C$.
 \medskip

\textbf{Claim~3.~} The only critical points of the function $\phi_0$ in the square $[-\frac{5}{8},\frac{5}{8}] \times [-\frac{\sqrt{3}}{3},\frac{\sqrt{3}}{2}]$
are the vertices $A, B, C$, and the mid-point $M_C$ of the side $AB$. \smallskip

We refer to \cite{BH-lmp16} for the explicit expression of the Neumann eigenvalues and eigenfunctions of the equilateral triangle $\cT_e$. After translation and rotation, we find that the second Neumann eigenfunction of $\cT_0$, which is symmetric with respect to $D_C$ is given by the formula,
\begin{equation}\label{E-appln-20}\scriptsize
\phi_0(u,v) = a_0\, \left( \cos\frac{4\pi u}{3} + \cos\frac{2\pi(1-u-\sqrt{3}v)}{3}
+ \cos\frac{2\pi(1+u-\sqrt{3}v)}{3} \right)\,,
\end{equation}
where $a_0\neq 0$ is a normalizing constant.

\begin{figure}
  \centering
  \includegraphics[scale=0.35]{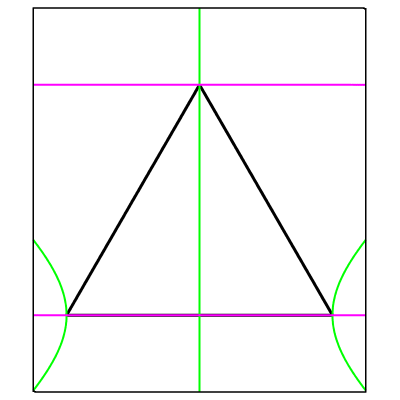}
  \caption{Localization of the critical points}\label{F-pthm-f6}
\end{figure}

It follows that the critical points of $\phi_0$ satisfy the equations,
\begin{equation}\label{E-appln-22}
\left\{
\begin{array}{l}
\sin\frac{2\pi u}{3} \, \left\lbrace  \cos\frac{2\pi u(1-\sqrt{3}v)}{3} + 2 \cos\frac{2\pi u}{3} \right\rbrace = 0\,,\\[5pt]
\sin\frac{2\pi u(1-\sqrt{3}v)}{3} \cos\frac{2\pi u}{3}  = 0\,.
\end{array}%
\right.
\end{equation}
The claim follows easily. It is also illustrated by Figure~\ref{F-pthm-f6} which displays the triangle $\cT_0$, the square $[-\frac{5}{8},\frac{5}{8}] \times [-\frac{\sqrt{3}}{3},\frac{\sqrt{3}}{2}]$, the zero set of $\partial_u\phi_0$ (green) and the zero set of $\partial_v\phi_0$ (magenta). \medskip

\textbf{Claim~4.~}  For $t$ small enough, $\phi_t+a < 0$ has exactly two connected components. \smallskip

For the proof, we proceed by contradiction. If not, there exists a sequence $t_n\rightarrow 0$, and a connected component $\omega(t_n)$ of $\phi_t+a < 0$, which according to Claim~2 must be contained in $\Omega(a,\varepsilon,\phi_0)$.\smallskip

Let $x_n \in \omega(t_n)$ be the point at which $\phi_{t_n}$ achieves its minimum in $\omega(t_n)$. We have necessarily $\nabla \phi_{t_n} (x_n)=0$. After extraction of a subsequence if necessary, we can assume that $x_n$ converged to some $x_\infty$ which belongs to $\wb{\cT_0}$, and satisfies $ -\epsilon \leq  \phi_0 (x_\infty) + a \leq \epsilon$. There are two possibilities. If $x_\infty \in \cT_0$, using Lemma~\ref{L-jnn-32s}, we get that $\phi_{t_n}$ converges to $\phi_0$ in a small ball around $x_\infty$ in $C^1$ sense, and this implies that $\nabla \phi_0 (x_\infty) =0$. A contradiction with Claim~3.\smallskip

The second possibility is that $x_\infty \in \partial \cT_0$. Here, we have to use a uniform boundary regularity for the Neumann Laplacian in $\Omega_t$ when we are far from $A,B,C$. We consider a small ball centered at
  $\partial \cT_0 \cap \{\phi_0 +a=0\}$ of radius $r(\epsilon)$ and containing  $\partial \cT_0 \cap \{-2\epsilon \leq \phi_0 +a \leq 2\epsilon \}$ (hence $x_\infty$). For each $t >0$, we consider a function $\chi (t,x)$ with support
   in the ball, equal to $1$ in a fixed neighborhood of $x_\infty$ and such that $\partial_\nu \chi (t,x) =0$ on $\partial \Omega_t$. It is easy to get such a function $C^\infty$ in both variables $t$ and $x$ due to the uniform regularity of $\partial \Omega(t)$ there (for $t\in [0,t_0]$ with $t_0>0$ small enough). We now consider $\hat \phi_t:=\chi(t,x) \phi_t$ in $\Omega_t$. This is a bounded family in $H^2$, and $\hat \phi_t$ satisfies the Neumann condition.

We have
   $$
   - \Delta \hat \phi_t = [-\Delta,\chi(t,\cdot)] \phi_t + \lambda_t \hat \phi_t\,.
   $$
The left hand side is uniformly bounded in $H^1$, and supported in the ball $B(x_\infty,r(\epsilon))$. We have a uniform (with respect to $t$) regularity of this Neumann problem (with locally  $C^\infty$ boundary), and we get that
the family $\hat \phi_t$ is bounded in $H^3(\Omega(t))$.

We now extend it in a bounded family  $\hat \Phi_t \in H^3_0(B(0,2M))$. Coming back to our sequence $\phi_{t_n}$, we observe that in particular $\hat \Phi_{t_n}$
   is a bounded family in $H^3_0(B(0,2M))$. Extracting a subsequence if necessary, we can assume that $\hat \Phi_{t_n}$ converges in $C^1(B(0,2M))$ to $\hat \Phi_\infty$. Now we have $\nabla \hat \phi_{t_n} (x_n)$ tends to $\nabla \hat  \Phi_\infty (x_\infty)$. For $n$ large enough $\nabla \hat \phi_{t_n} (x_n)=0$ which implies $\nabla \hat \Phi_\infty (x_\infty)=0$. Looking at the restriction to $\cT_0$, we also
    have $\hat \Phi_\infty =\chi (0,\cdot) \phi_0$ in $\cT_0$ in $\mathcal D'(\cT_0)$, which extends to $\wb{\cT_0}$ by continuity.\\
    This implies $0= \nabla \hat \Phi_\infty (x_\infty) = \nabla \phi_0(x_\infty)$, in contradiction with Claim~3. \hfill \qed\medskip

 \emph{Note.~} The preceding argument also shows that there cannot exist a second positive connected component for $t>0$ small enough  (without making use of the theorem of Gladwell and Zhu).

%\newpage
%\vspace{1cm}
\bibliographystyle{plain}

\end{document}